\documentclass[11pt,amsthm,reqno]{amsart}
\usepackage{amsmath,amssymb}

\newcommand{\C}{{\mathbb C}}
\newcommand{\Z}{{\mathbb Z}}

\newcommand{\A}{{\mathbb A}}

\newcommand{\Spec}{{\rm Spec}\;}
\newcommand{\Ker}{{\rm Ker}\:}

\newcommand{\codim}{{\rm codim}}
\newcommand{\Sing}{{\rm Sing}}
\newtheorem{thm}{Theorem}[section]
\newtheorem{cor}[thm]{Corollary}
\newtheorem{lem}[thm]{Lemma}

\newcommand{\Proof}{\noindent{\bf Proof.}\quad}
\newcommand{\Remark}{\noindent{\bf Remark.}\quad}

\begin{document}
\title[Factorial affine $G_a$-varieties]
{Factorial affine $G_a$-varieties with principal plinth ideals}
\author{Kayo Masuda}
\address{Department of Mathematical Sciences, 
School of Science, Kwansei Gakuin University \\
1 Gakuen Uegahara, Sanda 669-1330, Japan}
\email{kayo@kwansei.ac.jp}
\date{}

\begin{abstract}
Let $X=\Spec B$ be a factorial affine variety defined over an algebraically closed field $k$ of characteristic zero 
with a nontrivial action of the additive group $G_a$ associated to a locally nilpotent derivation $\delta$ on $B$.
Suppose that $A=\Ker \delta$ is an affine $k$-domain.
The quotient morphism $\pi : X \to Y=\Spec A$ splits to a composite ${\rm pr} \circ p$ 
of the projection ${\rm pr} : Y\times \A^1 \to Y$ and a $G_a$-equivariant birational morphism $p : X \to Y\times \A^1$ 
where $G_a$ acts on $\A^1$ by translation. 
In this article, we study $X$ of dimension $\ge 3$ under the assumption that 
the plinth ideal $\delta(B)\cap A$ is a principal ideal generated by a non-unit element $a$ of $A$. 
By decomposing $p : X \to Y\times \A^1$ to a sequence of $G_a$-equivariant affine modifications, we investigate the structure of $X$. 
We show in algebraic way that the general closed fiber of $\pi$ over the closed set $V(a)$ of $Y$ consists of a disjoint union of affine lines. 
The $G_a$-action on $X$ and the fixed-point locus $X^{G_a}$ are studied with particular interest. 
\end{abstract}

\maketitle

\section{Introduction} 

Let $k$ be an algebraically closed field of characteristic $0$, which is the ground field. 
An affine algebraic variety $X$ is factorial iff the algebra of regular functions on $X$ is a UFD. 
Note that the algebra of regular functions on a smooth acyclic affine variety is a UFD \cite{Fuj}.
Let $X=\Spec B$ be a factorial affine algebraic variety with a nontrivial algebraic action of the additive group $G_a$. 
There exists a nontrivial locally nilpotent derivation (abbreviated to lnd) $\delta$ on $B$ associated to the $G_a$-action on $X$. 
The fixed-point locus $X^{G_a}$ is defined by the ideal $(\delta(B))$ generated by the image $\delta(B)$. 
We assume that $\delta$ is irreducible, i.e., $\delta(B)$ is not contained in any proper principal ideal of $B$. 
If $\codim_X X^{G_a}>1$, then $\delta$ is irreducible. 
Let $A=\Ker \delta$. Then $A$ coincides with the invariant ring $B^{G_a}$. If $\dim X > 3$, $A$ is not necessarily finitely generated over $k$.
Suppose that $A$ is an affine $k$-domain and let $Y=\Spec A$.
The quotient morphism $\pi : X \to Y$ defined by the inclusion $A \hookrightarrow B$ is an $\A^1$-fibration
whose general closed fiber is a $G_a$-orbit $\A^1$. 
In fact, there exists an element $z \in B$ such that $B[a^{-1}]=A[a^{-1}][z]$
where $a=\delta(z)\in A\setminus \{0\}$ and $z$ is transcendental over $A[a^{-1}]$. 
Hence $\pi|_{\pi^{-1}(D(a))} : \pi^{-1}(D(a)) \to D(a)$ is a trivial $\A^1$-bundle
over the open set $D(a)=\{\mathfrak p\in \Spec A\mid \mathfrak p\not\ni a\}$. 
Let $p : X \to Y\times \A^1=\Spec A[z]$ be the $G_a$-equivariant morphism defined by the inclusion $A[z]\hookrightarrow B$.
Then $\pi$ splits to $\pi={\rm pr} \circ p$ where ${\rm pr} : Y\times \A^1 \to Y$ is the projection. 
In \cite{Fre}, \cite{Fre1}, Freudenburg defines the canonical factorization of the quotient morphism $\pi : X \to Y$,
which consists of ${\rm pr} : Y\times \A^1 \to Y$ and $G_a$-equivariant affine modifications developed by Kaliman and Zaidenberg \cite{KZ}. 
Inspired by the work of Kaliman-Zaidenberg and Freudenburg, the $G_a$-action on $X$ of dimension $\ge 3$ is studied in \cite{Ma1} 
by decomposing the morphism $p: X \to Y\times \A^1$ to a sequence of $G_a$-equivariant affine modifications, 
and given is a criterion for $X$ to be isomorphic to a hypersurface of type $x^m y-g(z)=0$ 
where $m>0$, $x$ is a prime element of $A$ and $g(z)\in A[z]\setminus A$. 
The ideal $\delta(B)\cap A$ of $A$ is called the {\it plinth ideal}. 
In this article, we study the structure of $X$ of dimension $\ge 3$ under the condition that the plinth ideal is principal. 
If $X$ and $Y$ are smooth and if $\pi : X \to Y$ is surjective and equi-dimensional, the plinth ideal is principal (Lemma \ref{Lemma 2.1}). 
Suppose that $\delta(B)\cap A=\alpha_1^{p_1}\cdots \alpha_s^{p_s}A$ 
where $\alpha_1, \ldots, \alpha_s$ are distinct prime elements of $A$ and $p_i>0$ for $1\le i \le s$. 
We decompose $p: X \to Y\times \A^1$ to a sequence of $G_a$-equivariant affine modifications 
$$X=X_s \to X_{s-1} \to \cdots \to X_1 \to Y\times \A^1.$$
By investigating the affine modifications, we show that the general closed fiber of $\pi : X \to Y$ over each $V(\alpha_i)=\Spec A/\alpha_i A$ 
consists of a disjoint union of $m_i$ affine lines for $m_i>1$ (Theorem \ref{Theorem 3.9}).
As a consequence, we obtain the following: Suppose that $X$ and $Y$ are smooth and
the quotient morphism $\pi : X\to Y$ is surjective and equi-dimensional.
Suppose, further, that the restriction $\pi |_{\pi^{-1}(D(a))}$ is a trivial $\A^1$-bundle over an open set $D(a)$ of $Y$. 
If the general closed fiber of $\pi$ over $V(a)=Y\setminus D(a)$ is irreducible, 
then $\pi : X \to Y$ is a trivial $\A^1$-bundle (Corollary \ref{Corollary 3.11}). 
Hence $X$ is $G_a$-equivariantly isomorphic to $Y \times \A^1$ where $G_a$ acts on $\A^1$ by translation. 
We also give a condition for $X$ to have no $G_a$-fixed points (Corollary \ref{Corollary 3.10}). 

For $n\ge 3$, an affine pseudo-$n$-space is, by definition, 
a smooth affine variety $Z$ equipped with a faithfully flat morphism $q : Z \to \A^1=\Spec k[x]$ 
such that $q^{-1}(\A^1_*)\cong \A^1_*\times \A^{n-1}$ and the scheme-theoretic fiber $q^*(0)$ is irreducible and reduced
where $\A^1_*=\A^1\setminus \{0\}$. 
If $X=\Spec B$ is an affine pseudo-$n$-space, then $X$ is factorial with $B^*=k^*$ and has a $G_a$-action associated to an irreducible lnd $\delta$
such that $k[x]\subset A=\Ker \delta$ and $\delta(B)\cap A\supset x^m A$ for some $m>0$ \cite{Ma1}.  
Suppose that $A$ is an affine $k$-domain and $\delta(B)\cap A$ is principal, and let $Y=\Spec A$. 
By applying Corollary \ref{Corollary 3.11} to the affine pseudo-$n$-space $X$, 
we obtain that $X$ is $G_a$-equivariantly isomorphic to $Y \times \A^1$ 
if the general closed fiber of the quotient morphism $\pi : X \to Y=\Spec A$ over $V(x)$ is irreducible (Corollary \ref{Corollary 3.13}). 
In particular, an affine pseudo-$3$-space $X=\Spec B$ such that $q^*(0)=\Spec B/xB$ is factorial and $(B/x B)^*=k^*$ 
is isomorphic to $Y \times \A^1\cong \A^3$ if the general closed fiber of $\pi$ over $V(x)$ is irreducible (Corollary \ref{Corollary 3.14}). 
Then $x$ is a variable of $X\cong \A^3$ by a result of Kaliman \cite{Kal2}. 
If the general closed fiber of $\pi$ over $V(x)$ is reducible,
then $X$ is not necessarily isomorphic to $Y\times \A^1$ although $X\cong \A^3$ (see Example 4.1). 
We illustrate the $G_a$-equivariant affine modifications with examples in the last section. 

\medskip

{\it Acknowledgement.} \; The author was supported by KAKENHI Grant Number JP20K03570 and JP20K03525, JSPS.

\medskip
\bigskip

\section{Preliminaries}

First, we recall some basic facts on lnds on an affine $k$-domain. 
We refer to Miyanishi \cite{Miy2} and Freudenburg \cite{Fre} for further details. 
For the quotient morphism by $G_a$-action, see also \cite{GMM3}. 

Let $B$ be an affine $k$-domain and let $\delta$ be a nontrivial lnd on $B$. 
Let $A=\Ker \delta$. The group $A^*$ of invertibles of $A$ coincides with $B^*$.
If $B$ is factorial, then $A$ is factorial as well and a prime element of $A$ is a prime element of $B$ 
since $A$ is factorially closed in $B$, i.e., 
$xy \in A\setminus \{0\}$ for $x,y \in B$ implies $x,y \in A$. 
There exists an element $z\in B$, called a local slice of $\delta$, which satisfies $\delta(z)=a\in A\setminus \{0\}$. 
It is well-known that $B[a^{-1}]=A[a^{-1}][z]$ and $z$ is transcendental over $A[a^{-1}]$.  
If $\delta(z)\in A^*$, $z$ is called a slice. 

An lnd $\delta$ is irreducible iff $\delta(B)\subset bB$ for $b \in B$ implies $b\in B^*$. 
There exist a nontrivial irreducible lnd $\bar \delta$ on $B$ and 
an element $a\in \Ker \bar\delta \setminus \{0\}$ 
such that $\delta=a\bar\delta$.  
Note that $\Ker \delta$ coincides with $\Ker \bar\delta$. 

Let $X=\Spec B$. The fixed point locus $X^{G_a}$ is defined by the ideal $(\delta(B))$
generated by the image $\delta(B)$. 
By a result of Bialynicki-Birula \cite{BB}, $X^{G_a}$ has no isolated fixed points. 
If $X$ is factorial and $\delta$ is irreducible, then $X^{G_a}$ has codimension $>1$. 

The ideal $\delta(B)\cap A$ of $A$ is called the plinth ideal. 
An ideal $I$ of $B$ is called {\it $\delta$-stable} or {\it integral} if $\delta(I)\subset I$. 

When $B=k^{[3]}$, $A=k^{[2]}$ by a result of Miyanishi \cite{Miy} and the plinth ideal $\delta(B)\cap A$ is principal
and the quotient morphism $\pi : X \to Y$ defined by the inclusion $A\hookrightarrow B$ is surjective
by Bonnet \cite{Bon} (in case $k=\C$), Daigle and Kaliman \cite{DK} (in case that $k$ is of characteristic $0$).

For elements $a_1,\ldots, a_r$ of a subdomain $R$ of $B$, we denote by $(a_1,\ldots,a_r)R$
(resp. $(a_1,\ldots, a_r)B$) the ideal of $R$ (resp. $B$) generated by $a_1,\ldots, a_r$.
We have the following result. 

\begin{lem}\label{Lemma 2.1}
Let $X=\Spec B$ be a smooth factorial affine variety with a $G_a$-action associated to an lnd $\delta$. 
Suppose that $A=\Ker \delta$ is an affine $k$-domain and $Y=\Spec A$ is smooth. 
If the quotient morphism $\pi : X \to Y$ is surjective and equi-dimensional,
then the plinth ideal $\delta(B)\cap A$ is principal.
\end{lem}
\Proof
Let $a_1, a_2$ be nonzero elements of $\delta(B)\cap A$ such that $a_i=\delta(s_i)$ where $s_i\in B$ for $i=1,2$.  
It suffices to show that $d=\gcd(a_1, a_2)\in \delta(B)\cap A$. 
Let $I=(a_1,a_2)A$. By the assumption, $B$ is faithfully flat over $A$. Hence we have $A\cap IB=I$. 
Since $\delta(a_1s_2-a_2s_1)=0$, it follows that $a_1s_2-a_2s_1 \in A\cap IB=I$. 
Hence $a_1s_2-a_2s_1=a_1c_1-a_2c_2$ for $c_1, c_2\in A$. 
Then we have $a_1(s_2-c_1)=a_2(s_1-c_2)$, hence $a'_1(s_2-c_1)=a'_2(s_1-c_2)$ where $a_i=da'_i$ for $i=1,2$. 
Define $s=(s_2-c_1)/a'_2=(s_1-c_2)/a'_1$. Then $s \in B$ satisfies $\delta(s)=d$.
Hence $d \in \delta(B)\cap A$, and the assertion follows.
\qed

\medskip

For $c\in B\setminus \{0\}$, we denote the localization of $B$ at $c$ by $B_c$ or $B[c^{-1}]$. 

\begin{lem}\label{Lemma 2.2}
Let $R$ be a subdomain of $B$ and $\alpha$ a nonzero element of $R$. 
\begin{enumerate}
\item[(1)] If $B[\alpha^{-1}]=R[\alpha^{-1}]$ and the ideal $R\cap \alpha B$ of $R$ is generated by $\alpha$, then $B=R$. 
\item[(2)] Suppose that $R\cap \alpha B=(\alpha, g_1, \ldots, g_r)R$ for $g_1, \ldots, g_r \in R$. 
If there exists some $\ell \ge 1$ such that $g_i\in R\cap \alpha^\ell B$ for $1 \le i \le r$, 
then $R\cap \alpha^j B=(\alpha^j, g_1, \ldots, g_r)R$ for $1 \le j \le \ell$. 
\end{enumerate} 
\end{lem}
\Proof
(1) Take any $b \in B\setminus \{0\}$. Then for a nonnegative integer $m$, $\alpha^m b=a\in R$. 
If $m>0$, then $a\in R\cap \alpha B=\alpha R$. Hence we have $\alpha^{m-1}b=a_1\in R$. 
If $m>1$, by repeating this argument, we have $b\in R$.

(2) Let $I_i=R\cap \alpha^i B$ for $i\ge 1$. We show $I_j=(\alpha^j, g_1,\ldots, g_r)R$ by induction on $j$ for $j\le \ell$.
Take any $h\in I_j$ for $1< j \le \ell$. Then since $h\in I_j \subset I_1$, $h$ is written as 
$h=\alpha h_0+g_1 h_1+\cdots +g_r h_r$ for $h_0, \ldots, h_r \in R$. 
Hence we have $\alpha h_0=h-(g_1 h_1+\cdots +g_rh_r)\in I_j$ and obtain
$h_0\in I_{j-1}=(\alpha^{j-1}, g_1, \ldots, g_r)R$.
Thus $h\in (\alpha^j, g_1, \ldots, g_r)R$, and the assertion follows. 
\qed

\medskip

Next, we review some basic facts of equivariant affine modifications of a factorial affine $G_a$-variety. 
We refer to Kaliman and Zaidenberg \cite{KZ} for details. 

Let $B$ be a factorial affine $k$-domain. 
Let $R$ be a factorial subdomain of $B$, $I$ a nontrivial ideal of $R$, and $f$ a nonzero element of $I$. 
The subalgebra of the quotient field $Q(R)$ generated over $R$ by the elements $a/f$ 
for $a\in I$ is denoted by $R[f^{-1}I]$ and called the {\it affine modification of $R$ along $f$ with center $I$}. 
If $I$ is generated by $a_1,\ldots, a_r$, then 
$$R[f^{-1}I]=R[a_1/f,\ldots, a_r/f].$$
In particular, if $I$ is generated by $f$ and $g$ which are coprime,
then $R[f^{-1}I]=R[g/f]$ is isomorphic to $R[Y]/(f Y-g)$ as an $R$-algebra where $R[Y]=R^{[1]}$ (cf. \cite{Ma1}). 
Further, if $R$ is noetherian and $I$ is a prime ideal, then $R[f^{-1}I]$ is factorial ({\it ibid.}). 

Suppose that $B$ is equipped with a nontrivial lnd $\delta$ restricting to $R$. 
If $I$ is $\delta$-stable and $\delta(f)=0$, then $\delta$ uniquely lifts up to the affine modification $R[f^{-1}I]$.
Suppose, further, that $B[f^{-1}]=R[f^{-1}]$. 
For $i\ge 1$, let $I_i$ be the ideal $R\cap f^i B$ of $R$ and let $B_i=R[f^{-i}I_i]$.
Then there exsits a finite sequence of $G_a$-equivariant affine modifications 
\begin{equation}
R=B_0\subset B_1\subset B_2 \subset \cdots \subset B_\mu=B. \label{(1)}
\end{equation}
where $\delta$ restricts to $B_i$ for $0 \le i \le \mu$ ({\it ibid}). 
If $B^{G_a}=R^{G_a}$, then $B_i^{G_a}=R^{G_a}$ for every $i$. 
The sequence (\ref{(1)}) of $G_a$-equivariant affine modifications yields a $G_a$-equivariant birational morphisms
$$X=X_\mu \to X_{\mu-1} \to \cdots \to X_1 \to X_0$$
where $X_i=\Spec B_i$ for $0 \le i \le \mu$. We also call $X_i$ a $G_a$-equivariant affine modification of $X_0$. 

Let $c\in R\setminus \{0\}$ be an element such that $\delta(c)=0$ and $f$ and $c$ are coprime.
By localizing the sequence (\ref{(1)}) at $c$, 
we have a sequence of $G_a$-equivariant affine modifications 
$$R_{c}=B_{0,c}\subset B_{1,c}\subset B_{2,c} \subset \cdots \subset B_{\mu, c}=B_c$$
where $B_{i,c}=B_{i}[c^{-1}]$ for $0 \le i \le \mu$. 
For each $i$, $B_{i,c}$ coincides with the $G_a$-equivariant affine modification of $R_c$ 
along $f^i$ with center $I_{i,c}=R_c\cap f^i B_c$. 

\begin{lem}\label{Lemma 2.3}
Let $B$ be a factorial affine $k$-domain with a nontrivial lnd $\delta$ which restricts to a factorial subdomain $R$ and  
$a\in R$ a nonzero element such that $\delta(a)=0$ and $B[a^{-1}]=R[a^{-1}]$. 
Write $a$ as $a=\alpha_1\alpha_2\cdots \alpha_s$ where $\alpha_1, \ldots, \alpha_s$ are pairwise coprime elements of $R$. 
For $1 \le i\le s$, let $I^{(i)}$ be the ideal $R\cap \alpha_1 \cdots \alpha_iB$ of $R$ and let $B^{(i)}=R[(\alpha_1\cdots\alpha_i)^{-1}I^{(i)}]$.
Then the following assertions hold. 
\begin{enumerate}
\item[(1)] There exsits a sequence of $G_a$-equivariant affine modifications 
\begin{equation}
R=B^{(0)}\subset B^{(1)}\subset B^{(2)} \subset \cdots \subset B^{(s)}\subset B. \label{(2)}
\end{equation}
where $\delta$ restricts to $B^{(i)}$ for $0 \le i \le s$. 
If $B^{G_a}=R^{G_a}$, then $(B^{(i)})^{G_a}=R^{G_a}$ for every $i$. 
\item[(2)] Let $\{b_1, \ldots, b_r\}$ be a set of generators of $B$ over $k$. 
If $a b_j \in R$ for $1 \le j \le r$, then $B^{(s)}=B$.
\end{enumerate}  
\end{lem}
\Proof
(1) We show $B^{(i)} \subset B^{(i+1)}$ for $1 \le i < s$. 
The affine modification $B^{(i)}$ is generated over $R$ by $c_l\in B$ such that $\alpha_1\cdots\alpha_i c_l=a_l \in R$.
Since $\alpha_1\cdots \alpha_i\alpha_{i+1}c_l=\alpha_{i+1}a_l\in R\cap \alpha_1\cdots\alpha_{i+1}B=I^{(i+1)}$, it follows that $c_l\in B^{(i+1)}$. 
Hence $B^{(i)} \subset B^{(i+1)}$.
Since $\delta(a)=\delta(\alpha_1\cdots \alpha_s)=0$, $\alpha_1\cdots\delta(\alpha_j)\cdots \alpha_s$ is a multiple of $\alpha_j$ for every $j$.
Then $\delta(\alpha_j)$ is a multiple of $\alpha_j$ since $\alpha_1, \ldots, \alpha_s$ are pairwise coprime.
Since $\delta$ is locally nilpotent, $\delta(\alpha_j)=0$ for every $j$ (cf. \cite{Fre}). 
Hence $I^{(i)}$ is $\delta$-stable and $\delta$ restricts to each $B^{(i)}$. 
If $B^{G_a}=R^{G_a}$, we have $(B^{(i)})^{G_a}=R^{G_a}$ by taking $G_a$-invariants of the sequence (\ref{(2)}). 

(2) Since $I^{(s)}=R\cap aB$ and $B^{(s)}=R[a^{-1}I^{(s)}]$, it follows from $ab_j\in R$ that $b_j \in B^{(s)}$ for every $j$. Hence $B^{(s)}=B$.
\qed

\medskip
\bigskip

\section{Equivariant affine modifications of factorial $G_a$-varieties}

Let $B$ be a factorial affine $k$-domain and let $\delta$ be a nontrivial irreducible lnd on $B$. 
Let $A=\Ker \delta$. Throughout this section, we assume that $A$ is noetherian and the plinth ideal is principal. 
Let
\begin{equation}
\delta(B)\cap A=\alpha_1^{p_1}\alpha_2^{p_2}\cdots \alpha_s^{p_s}A \label{(3)}
\end{equation}  
where $\alpha_i$ is a prime element of $A$, $\alpha_i\neq \alpha_j$ if $i \neq j$, and $p_i>0$ for $1\le i \le s$.  
Let $z \in B$ be a local slice such that  
$$\delta(z)=\alpha_1^{p_1}\alpha_2^{p_2}\cdots \alpha_s^{p_s}.$$ 
Then $z$ is not divisible by any $\alpha_i$. 
Note that $B\supsetneq A[z]$ since $\delta$ is irreducible. 

Let $b_1, \ldots, b_r$ be the generators of $B$ over $k$.
Since $B[a^{-1}]=A[a^{-1}][z]$ where $a=\alpha_1^{p_1}\alpha_2^{p_2}\cdots \alpha_s^{p_s}$,
there exist nonnegative integers $\mu_1, \ldots, \mu_s$ such that
$\alpha_1^{\mu_1}\alpha_2^{\mu_2}\cdots \alpha_s^{\mu_s}b_j\in A[z]$ for $1 \le j \le r$. 
We choose $\mu_i$ to be the minimal for every $i$. 

\begin{lem}\label{Lemma 3.1}
For $1 \le i \le s$, $\mu_i>0$ and $A[z]\cap \alpha_i B\supsetneq \alpha_i A[z]$. 
\end{lem}
\Proof
Suppose that $\mu_i=0$ for some $i$, say, $\mu_1=0$. 
Then $\alpha_2^{\mu_2}\cdots \alpha_s^{\mu_s}b_j=h_j(z) \in A[z]$ for $1 \le j \le r$. 
Hence for every $j$
$$\alpha_2^{\mu_2}\cdots \alpha_s^{\mu_s}\delta(b_j)=\alpha_1^{p_1}\alpha_2^{p_2}\cdots \alpha_s^{p_s}h'_j(z),$$
from which we have $\delta(b_j)$ is a multiple of $\alpha_1^{p_1}$.
This is a contradiction because $\delta$ is irreducible.
Hence $\mu_i>0$ for every $i$.

Suppose that $A[z]\cap \;\alpha_i B=\alpha_i A[z]$ for some $\alpha_i$, say $\alpha_1$.
Since $\alpha_1^{\mu_1}\cdots \alpha_s^{\mu_s}b_j \\ \in A[z]\cap \alpha_1 B=\alpha_1 A[z]$ for every $j$,
it follows that $\alpha_1^{\mu_1-1}\alpha_2^{\mu_2}\cdots \alpha_s^{\mu_s}b_j\in A[z]$,
which contradicts to the minimality of $\mu_1$.  
Hence the assertion follows. 
\qed

\medskip

For $1 \le i \le s$, let
$$I^{(i)}=A[z]\cap \alpha_1^{\mu_1}\cdots \alpha_i^{\mu_i}B \quad \text{and} \quad B^{(i)}=A[z][\alpha_1^{-\mu_1}\cdots \alpha_i^{-\mu_i}I^{(i)}].$$
Since $A[z]$ is a factorial subdomain of $B$, we have by Lemma \ref{Lemma 2.3} a sequence of $G_a$-equivariant affine modifications
\begin{equation}
A[z]\subset B^{(1)}\subset B^{(2)} \subset \cdots \subset B^{(s)}=B \label{(4)}
\end{equation}  
and $(B^{(i)})^{G_a}=A$ for every $i$. 

We investigate $B^{(1)}=A[z][\alpha_1^{-\mu_1}I^{(1)}]$ where $I^{(1)}=A[z]\cap \alpha_1^{\mu_1}B$.   
In the sequel, we denote $\alpha_1$ by $\alpha$, $p_1$ by $p$, and $\mu_1$ by $\mu$ for simplicity.
Let $\beta=\alpha_2^{\mu_2}\cdots \alpha_s^{\mu_s}$.
Then 
$$\delta(z)=\alpha^{p}\beta.$$
For $i\ge 1$, let 
$$I_{i}=A[z]\cap \alpha^{i}B \quad \text{and} \quad B_{i}=A[z][\alpha^{-i}I_{i}].$$
Then we have a sequence of $G_a$-equivariant affine modifications
\begin{equation}
A[z]\subset B_{1}\subset B_{2} \subset \cdots \subset B_{\mu}=B^{(1)} \label{(5)}
\end{equation}  
and $B_i^{G_a}=A$ for every $i$. 
We analize this sequence (\ref{(5)}) of $G_a$-equivariant affine modifications. 

By Lemma \ref{Lemma 3.1}, $I_1\supsetneq \alpha A[z]$. 
Note that $I_1$ is the prime ideal of $A[z]$.
Further, $I_1 \cap A=\alpha A$ since $A$ is factorially closed in $B$.

Let $\overline A=A/\alpha A$. 
The residue ring $A[z]/\alpha A[z]$ is identified with a polynomial ring $\overline A[\overline z]$ over $\overline A$ 
where $\overline z$ is the residue class of $z$. 
For an ideal $I$ of $A[z]$, we denote by $\overline I$ the image of $I$ 
by the surjection $A[z] \to A[z]/\alpha A[z]=\overline A[\overline z]$. 
Then $\overline I_1\neq(0)$, $\overline I_1\cap \overline A=(0)$ and 
$\overline I_1$ is a prime ideal of $\overline A[\overline z]$. 

Let $K$ be the quotient field $Q(\overline A)$ of $\overline A$.
For $i\ge 1$, let ${\overline I}_i^K=\overline I_i\otimes_{\overline A}K$.
The ideal $\overline I_j^K$ of $K[\overline z]$ is principal and
satisfies $\overline I_i^K \supset \overline I_j^K$ for $i\le j$.
Suppose that for some $1 \le \ell_1 \le \mu$
\begin{equation*}
{\overline I}_1^K = \overline I_{\ell_1}^K\supsetneq \overline I_{\ell_1+1}^K\supset \cdots \supset \overline I_\mu^K.
\end{equation*}
Let $g\in I_{\ell_1}\subset I_1$ be an element which maps to a generator of $\overline I_{\ell_1}^K=\overline I_1^K$
by the map $A[z] \to \overline A[\overline z]\hookrightarrow K[\overline z]$. 
Note that $g\notin A$ since $\overline I_1\neq (0)$. 
Since $\overline I_1^K$ is a prime ideal, $\overline g\in \overline I_{\ell_1}\subset \overline I_1$ is
an irreducible polynomial in $K[\overline z]$.
We may assume that $g=g(z)\in A[z]$ is primitive over $A$. 
Since $\overline g$ is a generator of $\overline I_{\ell_1}^K=\overline I_1^K$, 
$\overline g\in \overline A[\overline z]$ has the minimal degree in $\overline I_1\setminus \{0\}$ with respect to $\overline z$.
Hence $g=g(z)\in A[z]$ is irreducible since $g(z)$ is contained in the prime ideal $I_1$ and primitive over $A$. 
Write $g(z)\in A[z]\setminus A$ as
$$g(z)=\tilde g(z)+\alpha h(z)$$
where $\tilde g(z)\in I_1\setminus \alpha A[z]$ and $h(z) \in A[z]$.
Let $c\in A$ be the coefficient of the highest term of $\tilde g(z)$. We may assume $c\notin \alpha A$. 
Then $\tilde g(z)$ has the minimal degree with respect to $z$ in $I_1\setminus \alpha A[z]$
since $\overline{\tilde g}=\overline g\in \overline I^K_1$. 
By localizing the sequence (\ref{(5)}) at $c\in A\setminus \alpha A$, we have a sequence of $G_a$-equivariant affine modifications
\begin{equation*}
A_c[z]\subset B_{1,c}\subset \cdots \subset B_{\mu,c}=B^{(1)}_c 
\end{equation*}  
where $B_{i,c}=B_i[c^{-1}]=A_c[z][\alpha^{-i}I_{i,c}]$ with $I_{i,c}=A_c[z]\cap \alpha^iB_c$ for $1 \le i \le \mu$. 
We also write $g(z)\in I_{\ell_1}$ as 
\begin{equation}
\qquad g(z)=\alpha^{\ell_1} y_1 \label{(6)}  
\end{equation} 
for $y_1\in B$. Note that $y_1 \notin \alpha B$ since $\overline g \notin \overline I^K_{\ell_1+1}$. 

\begin{lem}\label{Lemma 3.2}
With the notation above, the following assertions hold. 
\begin{enumerate}
\item[(1)] The ideal $I_{1,c}$ of $A_c[z]$ is generated by $\alpha$ and $g \in A[z]\setminus \alpha A[z]$. 
If $\overline A$ is factorial, then $I_1=(\alpha, g)A[z]$. 
\item[(2)] $\deg_{\overline z}\overline g(\overline z)>1$. 
\item[(3)] For any $a\in A$, $g'(z)-a \notin I_1$. 
\item[(4)] $q_1:=p-\ell_1\ge 0$ and $\delta(y_1)=\alpha^{q_1}\beta g'(z)$. 
\item[(5)] $I_{\ell_1,c}=(\alpha^{\ell_1}, g)A_c[z]$ and 
$$B_{\ell_1,c}=A_{c}[z,y_1]\cong A_c[z][Y]/(\alpha^{\ell_1} Y-g)$$ 
where $Y$ is an indeterminant.
If $s=1$, i.e., $\delta(z)=\alpha^p$, then $B_{\ell_1,c}$ is factorial and $\alpha$ is a prime element of $B_{\ell_1,c}$. 
\end{enumerate}
\end{lem}
\Proof
(1) Every $h(z)\in I_{1,c}$ is written as $h=\tilde g q+r$ where $q, r \in A_c[z]$ and $\deg_z r<\deg_z \tilde g$. 
Since $r=h-\tilde g q\in I_{1,c}$, it follows that $r\in \alpha A_c[z]$ 
by the minimality of the degree of $\tilde g$ in $I_{1,c}\setminus \alpha A_c[z]$. 
Hence $h\in (\alpha, g)A_c[z]$ and the assertion follows.
If $\overline A$ is factorial, the assertion follows from \cite[Lemma 4.1]{Ma1}.

(2) Suppose that $\deg \overline g(\overline z)=1$. Then $\tilde g(z)\in I_1$ is written as $\tilde g(z)=cz+c_0$ for $c_0 \in A$.
Since $\tilde g(z)=\alpha b$ for a nonzero $b\in B$, we have $\delta(b)=c\alpha^{p-1}\beta \in \delta(B)\cap A$. 
This contradicts to the equation (\ref{(3)}). 

(3) Since $\deg_z \tilde g'\ge 1$ by (2), $g'-a={\tilde g}'+\alpha h'-a\notin \alpha A[z]$ and $\deg_z (g'-a) <\deg_z g$. 
The assertion follows from that $\overline g$ has the minimal degree in $\overline I_1\setminus \{0\}$. 

(4) We show $\ell_1\le p$. Suppose the contrary. Applying $\delta$ to the equation (\ref{(6)}), we have
$$\beta g'=\alpha^{\ell_1-p}\delta(y_1)\in A[z]\cap \alpha B=I_1.$$
Then it follows that $g'\in I_1$, which is a contradiction by (3). 
Hence we have $p\ge \ell_1$ and the expression of $\delta(y_1)$. 

(5) The first assertion follows from (1) and Lemma \ref{Lemma 2.2}(2). 
Since $\alpha$ and $g$ are coprime in $A_c[z]$,
$B_{\ell_1,c}=A_c[z,y_1]$ is isomorphic to $A_c[z][Y]/(\alpha^{\ell_1} Y-g)$ as an $A_c[z]$-algebra (\cite{Ma1} cf. \cite{Ma}).
If $\delta(z)=\alpha^p$, we have $B_{\ell_1,c}[\alpha^{-1}]=A_{c}[\alpha^{-1}][z]=B_{c}[\alpha^{-1}]$. 
Hence $B_{\ell_1,c}[\alpha^{-1}]$ is factorial. 
We have
\begin{eqnarray*}
B_{\ell_1,c}/\alpha B_{\ell_1,c} &=& A_{c}[z,y_1]/\alpha A_{c}[z,y_1] \\
 &\cong & (A_{c}[z]/(\alpha, g)A_{c}[z])[\overline y_1] \\ 
&= & (A_{c}[z]/I_{1,c})[\overline y_1].
\end{eqnarray*}
Since $I_{1,c}=A_{c}[z]\cap \alpha B_{c}$ is a prime ideal of $A_{c}[z]$, $B_{\ell_1,c}/\alpha B_{\ell_1,c}$ is an integral domain.
Hence $\alpha$ is a prime element of $B_{\ell_1,c}$. 
Thus $B_{\ell_1,c}$ is factorial by a result of Nagata \cite{Nag}. 
\qed

\medskip

By Lemma \ref{Lemma 3.2}, We have
$$A_c[z]\subsetneq B_{\ell_1,c}=A_c[z,y_1]\subset B_{\mu,c}=B^{(1)}_c.$$
Suppose  
\begin{align*}
\overline I_{\ell_1}^K=\overline gK[\overline z]\supsetneq \overline I_{\ell_1+1}^K=\cdots = & \overline I_{\ell_2}^K
\supsetneq \overline I_{\ell_2+1}^K= \cdots =\overline I_{\ell_3}^K\supsetneq \overline I_{\ell_3+1}^K \cdots \notag \\
\cdots & \overline I_{\ell_{m-1}}^K \supsetneq \overline I_{\ell_{m-1}+1}^K =\cdots
=\overline I_{\ell_m}^K=\overline I_\mu^K.
\end{align*}
For $\ell_1\le j\le \ell_m$, let $g_{j}\in I_{j}$ be an element which maps to a generator of $\overline I_{j}^K$
by the map $A[z]\to \overline A[\overline z]\hookrightarrow K[\overline z]$. 
We take $g_{\ell_1}=g$. Since $g_\mu\in I_\mu\subset I_{\ell_m}$, we take $g_j=g_\mu$ for $\ell_m \le j \le \mu$. 
If $\overline I_j^K=\overline g^{e_j}K[\overline z]$ for some $e_j>0$ and $j\le e_j\ell_1$, we can take $g_{j}=g^{e_j}$. 

\begin{lem}\label{Lemma 3.3}
For $2\le j\le m$, $\overline I_{\ell_j}^K=\overline g^{j}K[\overline z]$ and $j\ell_1\le \ell_j$. 
\end{lem}  
\Proof
We first show that $\overline I_{\ell_j}^K=\overline g^{e_j}K[\overline z]$ for a positive integer $e_j$. 
For $2\le j\le m$, let $\xi_j\in K[\overline z]$ be a generator of $\overline I_{\ell_j}^K$. 
Since $\overline I_{\ell_1}^K \supset \overline I_{\ell_j}^K$, $\xi_j$ is written as
$\xi_j=\overline g \eta_j$ for $\eta_j\in K[\overline z]$.
While, since $g^e\in I_{\ell_j}$ for a sufficiently large $e$,
we have $\overline g^e=\xi_j \theta_j$ for $\theta_j \in K[\overline z]$.
Hence we have $\overline g^e=\overline g \eta_j \theta_j$.
Since $\overline g$ is irreducible in $K[\overline z]$,
$\eta_j$ and $\theta_j$ are some powers of $\overline g$ up to units. 
Hence $\overline I_{\ell_j}^K=\overline g^{e_j}K[\overline z]$ for some $e_j\ge 1$.
Note that $e_j<e_{j+1}$ since $\overline I_{\ell_j}^K\supsetneq \overline I_{\ell_{j+1}}^K$. 

Since $2\ell_1\ge \ell_1+1$, $g^2\in I_{\ell_1+1}$.
It follows from $\overline g^2\in \overline I_{\ell_1+1}^K=\overline I_{\ell_2}^K=\overline g^{e_2}K[\overline z]$
that $e_2=2$. Then $e_j=j$ holds by the induction on $j$.
In fact, suppose $e_j=j$. Then $\overline g_{\ell_j}=\gamma\overline g^{j}$ in $K[\overline z]$ for $\gamma \in K^*$.
Since $g_{\ell_j}g\in I_{\ell_j+1}$, it follows that 
$\overline g_{\ell_j}\overline g=\gamma\overline g^{j+1}\in \overline I^K_{\ell_j+1}=\overline I^K_{\ell_{j+1}}
=\overline g^{e_{j+1}}K[\overline z]$. Hence $j=e_j<e_{j+1} \le j+1$, and $e_{j+1}=j+1$.
The assertion $j\ell_1\le \ell_j$ follows from $\overline g^j\in \overline I^K_{j\ell_1}$ and
$\overline I^K_{\ell_j}=\overline g^jK[\overline z]\supsetneq \overline I^K_{\ell_j+1}=\overline I^K_{\ell_{j+1}}=\overline g^{j+1}K[\overline z]$.  
\qed

\medskip

By Lemmas \ref{Lemma 2.2} and \ref{Lemma 3.2},  
$$I_{\ell_1,c}=(\alpha^{\ell_1}, g)A_c[z]=\alpha I_{\ell_1-1,c}+ gA_c[z].$$
We set $I_{0,c}=A_c[z]$. 

\begin{lem}\label{Lemma 3.4}
There exists $c\in A\setminus \alpha A$ such that 
\begin{enumerate}
\item[(1)] $I_{j,c}=\alpha I_{j-1,c}+g_j A_c[z]$ for $\ell_1 \le j \le \mu$,
\item[(2)] $g_{\ell_j}=c_jg^{j}+\alpha f_j$ for $1<j \le m$ 
where $f_j \in I_{j\ell_1-1,c}$ and $c_j\in A_c^*$.  
\end{enumerate}
\end{lem}
\Proof
We show that there exists $d_j\in A\setminus \alpha A$ such that
$I_{j,d_j}=A_{d_j}[z]\cap \alpha^j B_{d_j}\subset \alpha I_{j-1,d_j}+g_{j}A_{d_j}[z]$.
For $j=\ell_1$, we take $d_{\ell_1}=c$, the coefficient of the highest term of $\tilde g(z)\in A[z]$. 
Let $j>\ell_1$ and let $h_1, \ldots, h_l$ be the generators of $I_{j}$. 
Since $\overline I_{j}^K=\overline g_j K[\overline z]$, we have 
$\overline h_i=\gamma_i \overline g_j\overline q_i$ in $K[\overline z]$ for $1 \le i \le l$
where $\gamma_i\in K^*$ and $q_i\in A[z]$. 
Thus $a_ih_i=a'_ig_{j}q_i+\alpha r_i$ for $a_i, a_i'\in A\setminus \alpha A$ and $r_i\in A[z]$.
Since $h_i, g_j\in I_{j}$, it follows that $\alpha r_i\in I_{j}$. 
Hence $r_i\in I_{j-1}$, and $a_i h_i\in \alpha I_{j-1}+g_{j} A[z]$. 
Let $d_j=a_1\cdots a_l$. Then we have $I_{j,d_j}\subset \alpha I_{j-1,d_j}+g_{j} A_{d_j}[z]$. 
By setting $c=\prod_{j=\ell_1}^{\mu}d_j$, it holds that $I_{j,c}\subset \alpha I_{j-1,c}+g_j A_c[z]$ for every $j$, 
and hence (1) is satisfied. 

Since $\overline g_{\ell_j}=\beta_j \overline g^{j}$ for $\beta_j \in K^*$,
we have $b_j g_{\ell_j}=b_j'g^{j}+\alpha f_j$ where $b_j, b_j'\in A\setminus \alpha A$ and $f_j\in A[z]$. 
Then $f_j \in I_{j\ell_1-1}$ since $\alpha f_j=b_j g_{\ell_j}-b_j'g^{j}\in I_{j\ell_1}$. 
Replacing $c$ by $c\prod_{j=2}^m b_jb_j'$, (2) is also satisfied, and we obtain a required $c\in A\setminus \alpha A$. 
\qed

\medskip

In the sequel, $c$ denotes an element of $A\setminus \alpha A$ satisfying the conditions (1) and (2) in Lemma \ref{Lemma 3.4}. 
By Lemma \ref{Lemma 3.4}, $\overline g_{\ell_j}=\overline c_j\overline g^j$ for $2 \le j \le m$ where $c_j\in A_c^*$.
We have a sequence of ideals of $A_c[z]$
$$I_{\ell_1,c}\supset I_{\ell_2,c}\supset \cdots \supset I_{\ell_m,c} \supset I_{\mu,c}.$$
By Lemma \ref{Lemma 3.4}, $I_{\mu,c}=\alpha^{\mu-\ell_m}I_{\ell_m,c}+g_{\ell_m}A_c[z]$ since $g_j=g_\mu$ for $\ell_m \le j \le \mu$. 
Hence $B^{(1)}_c=B_{\mu,c}=B_{\ell_m,c}$.
There is a sequence of $G_a$-equivariant affine modifications
\begin{equation}
A_c[z]\subsetneq B_{\ell_1,c}=A_c[z,y_1]\subset B_{\ell_{2},c}\subset \cdots \subset B_{\ell_m,c}=B^{(1)}_c. \label{(7)}
\end{equation}

\begin{lem}\label{Lemma 3.5}
For $2 \le j \le m$, the following assertions hold. 
\begin{enumerate}
\item[(1)] For $\ell_{j-1} \le i \le \ell_j$, $I_{i,c}=\alpha^{i-\ell_{j-1}}I_{\ell_{j-1},c}+g_{\ell_j}A_c[z]$. 
\item[(2)] Suppose $\ell_{j-1} \le j\ell_1 \le \ell_j$. 
For $\ell_{j-1}\le i \le j\ell_1$, $I_{i,c}=\alpha^{i-\ell_{j-1}} I_{\ell_{j-1},c}+g^{j} A_c[z]$. 
For $j\ell_1 \le i \le \ell_j$, $I_{i,c}=\alpha^{i-j\ell_1}I_{j\ell_1,c}+g_{\ell_j}A_c[z]$. 
\item[(3)] Suppose that there exists $t\ge 1$ such that $\ell_s=s\ell_1$ for every $s$ such that $t < s \le j$. Then
\begin{align*}
I_{\ell_j,c} & = \alpha^{j\ell_1-\ell_t}I_{\ell_t,c} \\
& \qquad +(\alpha^{(j-t-1)\ell_1}g^{t+1}, \alpha^{(j-t-2)\ell_1}g^{t+2}, \ldots, \alpha^{\ell_1}g^{j-1}, g^j)A_c[z]
\end{align*}
and $B_{\ell_t,c}=B_{\ell_s,c}$ for $t \le s \le j$.   
In particular, if $\ell_s=s\ell_1$ for every $s$ such that $2 \le s \le j$, then 
$$I_{\ell_j,c}=(\alpha^{j\ell_1}, \alpha^{(j-1)\ell_1}g, \ldots, \alpha^{\ell_1}g^{j-1}, g^{j})A_c[z]$$
and $B_{\ell_1,c}=B_{\ell_2,c}=\cdots=B_{\ell_{j},c}=B_{(j+1)\ell_1,c}$. 
\item[(4)] If $\ell_j = e_1\ell_{j_1}+e_2\ell_{j_2}+\cdots +e_r\ell_{j_r}$ 
for some positive integers $e_1, \ldots, e_r$ and $j_1,\ldots, j_r$ such that $e_1 j_1+\cdots +e_r j_r=j$, then
$$I_{\ell_j,c}=\alpha^{\ell_j-\ell_{j-1}} I_{\ell_{j-1},c} \\ +g_{\ell_{j_1}}^{e_1}\cdots g_{\ell_{j_r}}^{e_r}A_c[z].$$
Hence $B_{\ell_{j-1},c}=B_{\ell_j,c}$. 
\end{enumerate}
\end{lem}
\Proof
(1) Since we can take $g_i=g_{\ell_j}$ for $\ell_{j-1}<i\le \ell_j$, the assertion follows from Lemma \ref{Lemma 3.4} (1).
The equation holds for $i=\ell_{j-1}$ as well since $g_{\ell_j}\in I_{\ell_{j-1},c}$. 

(2) We can take $g_i=g^{j}$ for $\ell_{j-1}<i \le j\ell_1$ and $g_i=g_{\ell_j}$ for $j\ell_1<i \le \ell_j$. 
The assertion follows from Lemma \ref{Lemma 3.4} (1). 
Note that the first assertion holds for $i=\ell_{j-1}$ as well since $g^j\in I_{\ell_{j-1},c}$.
Similarly, the second assertion holds for $i=j\ell_1$ as well since $g_{\ell_j}\in I_{j\ell_1,c}$.

(3) The first assertion follows by using (1) iteratively. 
We have by (2)
$$I_{(j+1)\ell_1,c}=\alpha^{\ell_1}I_{\ell_j,c}+g^{j+1}A_c[z]
=(\alpha^{(j+1)\ell_1}, \alpha^{j\ell_1}g, \ldots, \alpha^{\ell_1}g^{j}, g^{j+1})A_c[z].$$ 
Hence $B_{\ell_1,c}=\cdots=B_{\ell_j,c}=B_{(j+1)\ell_1,c}$.  

(4) By the assumption, we can take $g_{\ell_j}=g_{\ell_{j_1}}^{e_1}\cdots g_{\ell_{j_r}}^{e_r}$
since $\overline g_{\ell_t}=\overline c_t \overline g^t$ for $t=1, \cdots, r$. 
By (1), the assertion follows. 
\qed

\medskip

Let $t_1=1$ and let $t_2$ be a positive integer such that $\ell_s=s\ell_1$ for every $s$ such that $t_1 < s < t_2$ and $\ell_{t_2}>t_2\ell_1$. 
Then for any $s$ such that $t_1\le s <t_2$, $B_{\ell_1,c}=B_{\ell_s,c}\subset B_{\ell_{t_2},c}$ by Lemma \ref{Lemma 3.5} (3). 
For $j\ge 3$, we inductively define $t_j$ to be a positive integer such that
\begin{enumerate}
\item[(1)] for any $t_{j-1}< s <t_j$, there exsit nonnegative integers $e_{s_1}, \ldots, e_{s_{j-1}}$ 
satisfying $e_{s_1}+e_{s_2}t_2+\cdots +e_{s_{j-1}}t_{j-1}=s$ and $\ell_s=e_{s_1}\ell_1+e_{s_2}\ell_{t_2}+\cdots +e_{s_{j-1}}\ell_{t_{j-1}}$,
\item[(2)] $\ell_{t_j}>e_1\ell_1+e_2\ell_{t_2}+\cdots +e_{j-1}\ell_{t_{j-1}}$ holds for any nonnegative integers $e_1, \ldots, e_{j-1}$ satisfying 
$e_1+e_2t_2+\cdots +e_{j-1}t_{j-1}=t_j$.
\end{enumerate}  
Then $B_{\ell_{t_{j-1}},c}=B_{\ell_{t_j-1},c}\subset B_{\ell_{t_j},c}$ by Lemma \ref{Lemma 3.5} (4) and 
we obtain a subsequence
$$\ell_1=\ell_{t_1}<\ell_{t_2}< \cdots <\ell_{t_\nu}$$
of $\ell_1<\ell_2<\cdots <\ell_m$. 
Write $g_{\ell_{t_j}}\in I_{\ell_{t_j}}$ as 
\begin{equation}
g_{\ell_{t_j}}=\alpha^{\ell_{t_j}}y_j \label{(8)}
\end{equation}  
where $y_j\in B$. Note that $y_j\notin \alpha B$ since $g_{\ell_{t_j}}\in I_{\ell_{t_j},c}\setminus I_{\ell_{t_j}+1,c}$. 
By Lemma \ref{Lemma 3.5} (1),  we have
$$B_{\ell_{t_j},c}=B_{\ell_{t_{j-1}},c}[y_j]=A_c[z,y_1,\ldots, y_j]$$
and obtain a subsequence of (\ref{(7)}) 
$$A_c[z]\subsetneq B_{\ell_1,c}\subset B_{\ell_{t_2},c}\subset B_{\ell_{t_3},c}\subset \cdots \subset B_{\ell_{t_\nu},c}=B^{(1)}_c.$$

In the sequel, we assume $B_{\ell_1,c}\subsetneq B^{(1)}_c$. For $B_{\ell_{t_2},c}$, we have the following. 

\begin{lem}\label{Lemma 3.6}
\begin{enumerate}
\item[(1)] Let $\tilde q_2=\ell_{t_2}-t_2\ell_1$. Then 
\begin{equation} 
\alpha^{\tilde q_2}y_2=\tilde h_2(z,y_1)\in B_{\ell_1,c}=A_c[z,y_1] \label{(9)}
\end{equation}  
for $\tilde h_2(z,y_1)=c_{t_2}y_1^{t_2}+r_{t_2-1}(z)y_1^{t_2-1}+\cdots +r_1(z)y_1+r_0(z)+\alpha a_2(z,y_1)$ 
where $a_2(z,y_1)\in B_{\ell_1,c}$ and $r_j(z)\in A_c[z]\setminus \alpha B_c$ unless $r_j(z)=0$ for $0 \le j \le t_2-1$. 
\item[(2)] It holds that $B_{\ell_1,c}\cap \alpha B_c=(\alpha, \tilde h_2)B_{\ell_1,c}$. 
Let $q_2=q_1-\tilde q_2$. Then $q_2\ge 0$ and 
$$\delta(y_2)=\alpha^{q_2}\beta g' h_2\quad \text{modulo} \; \alpha^{q_2+1}\beta B_{\ell_1,c}$$ 
where $h_2=t_2c_{t_2}y_1^{t_2-1}+(t_2-1)r_{t_2-1}(z)y_1^{t_2-2}+\cdots +r_1(z)\in B_{\ell_1,c}\setminus \alpha B_c$. 
\end{enumerate}  
\end{lem}
\Proof
(1) By Lemma \ref{Lemma 3.4}, we have 
$g_{\ell_{t_2}}=c_{t_2}g^{t_2}+\alpha f_{t_2}$ for $c_{t_2}\in A_c^*$ and $f_{t_2}\in I_{t_2\ell_1-1,c}$. 
Since $\ell_{t_2-1}=(t_2-1)\ell_1$, $\ell_{t_2-1}< t_2\ell_1 < \ell_{t_2}$. 
Then $I_{t_2\ell_1-1,c}=\alpha^{\ell_1-1}I_{\ell_{t_2-1},c}+g^{t_2}A_c[z]$ by Lemma \ref{Lemma 3.5} (2). 
Hence by Lemma \ref{Lemma 3.5} (3), 
$$I_{t_2\ell_1-1,c}=(\alpha^{t_2\ell_1-1}, \alpha^{(t_2-1)\ell_1-1}g, \ldots, \alpha^{\ell_1-1}g^{t_2-1}, g^{t_2})A_c[z].$$
Since $g=\alpha^{\ell_1}y_1$, $g_{\ell_{t_2}}$ is written as 
$$\alpha^{\ell_{t_2}}y_2=g_{\ell_{t_2}}=c_{t_2}g^{t_2}+\alpha f_{t_2}=\alpha^{t_2\ell_1}a_1(z,y_1)+\alpha^{t_2\ell_1+1}y_1^{t_2}r(z)$$
where $a_1(z,y_1)=c_{t_2}y_1^{t_2}+r_{t_2-1}(z)y_1^{t_2-1}+\cdots +r_1(z)y_1+r_0(z)$ 
with $r(z), r_0(z),\\ \ldots, r_{t_2-1}(z) \in A_c[z]$. 
Hence we have
\begin{equation*}
\alpha^{\ell_{t_2}-t_2\ell_1}y_2=a_1(z,y_1)+\alpha a_2(z,y_1)\in B_{\ell_1,c} 
\end{equation*}
where $a_2(z,y_1)\in B_{\ell_1,c}$. 
We may assume $r_j(z) \notin \alpha B_c$ unless $r_j(z)=0$ for $0 \le j \le t_2-1$, and the assertion follows. 

(2) By (1), it follows that $a_1(z,y_1)\in B_{\ell_1,c}\cap \alpha B_c$. 
We show $B_{\ell_1,c}\cap \alpha B_c=(\alpha, a_1)B_{\ell_1,c}$. 
Let $u(z,y_1)$ be a nonzero element of $B_{\ell_1,c}\cap \alpha B_c$. 
Since $a_1=c_{t_2}y_1^{t_2}+\text{(terms with degree $\le t_2-1$ w.r.t. $y_1$)}$ with $c_{t_2}\in A_c^*$, $u(z,y_1)\in A_c[z,y_1]$ is written as
$$u(z,y_1)=b_0(z)+b_1(z)y_1+\cdots +b_{t_2-1}(z)y_1^{t_2-1}+\alpha q_0(z,y_1)+ a_1 q_1(z,y_1)$$
where $q_0(z,y_1), q_1(z,y_1)\in B_{\ell_1,c}$ and $b_0(z), \ldots, b_{t_2-1}(z)\in A_c[z]\setminus \alpha B_c$ unless zero. 
We have $b_0(z)+\cdots +b_{t_2-1}(z)y_1^{t_2-1}\in \alpha B_c$ since $a_1, u(z,y_1)\in \alpha B_c$. 
Suppose that $b_j(z)\neq 0$ and $b_{j+1}(z)=\cdots =b_{t_2-1}(z)=0$ for some $j$. 
Then $j>0$ and
$\alpha^{j\ell_1}b_0(z)+\cdots +\alpha^{\ell_1}b_{j-1}(z)g^{j-1}+b_j(z)g^j=\alpha^{j\ell_1}(b_0(z)+\cdots +b_j(z)y_1^j)\in I_{j\ell_1+1,c}$. 
Since $j<t_2$, $\ell_j=j\ell_1$. Hence it follows that
$$\overline b_j\overline g^j \in \overline I_{j\ell_1+1}^K=\overline I_{\ell_j+1}^K=\overline I_{\ell_{j+1}}^K=\overline g^{j+1}K[\overline z],$$ 
which is a contradiction since $\overline b_j\notin \overline g \overline A[\overline z]$. 
Thus we have $u(z,y_1)=\alpha q_0(z,y_1) \\+a_1 q_1(z,y_1)$,
and $B_{\ell_1,c}\cap \alpha B_c=(\alpha, a_1)B_{\ell_1,c}=(\alpha, \tilde h_2)B_{\ell_1,c}$. 

By (\ref{(9)}), we have 
\begin{align*}
\quad & \alpha^{\tilde q_2}\delta(y_2) \\
= & (\partial_{y_1} a_1(z,y_1)+\alpha \partial_{y_1} a_2(z,y_1))\delta(y_1)
+(\partial_{z} a_1(z,y_1)+\alpha \partial_{z} a_2(z,y_1))\delta(z) \\
= & \alpha^{q_1}\beta g'\partial_{y_1}a_1+\alpha^{q_1+1}\beta \xi_2
\end{align*}
where $\xi_2 \in B_{\ell_1,c}$. 
By the argument above, $\partial_{y_1}a_1 \notin \alpha B_c$ since $\deg_{y_1}\partial_{y_1}a_1<t_2$. 
Suppose $\tilde q_2 >q_1$. 
Then we have $g'\partial_{y_1}a_1 \in \alpha B_c$, which is a contradiction since $g', \partial_{y_1}a_1 \notin \alpha B_c$.
Hence $\tilde q_2\le q_1$, and we obtain the expression of $\delta(y_2)$. 
\qed

\medskip

By Lemma \ref{Lemma 3.6},
$$\delta(g_{\ell_{t_2}})=\alpha^{\ell_{t_2}}\delta(y_2)=\alpha^{p+(t_2-1)\ell_1}\beta(g'h_2+\alpha b_1)$$
for $b_1\in B_{\ell_1,c}$. 
While, $\delta(g_{\ell_{t_2}})=g'_{\ell_{t_2}}\delta(z)=\alpha^p\beta g'_{\ell_{t_2}}$.
Hence we have
$$g'_{\ell_{t_2}}=\alpha^{(t_2-1)\ell_1}(g'h_2+\alpha b_1).$$ 

Let $3 \le j \le \nu$. Let $\sigma_{j,j-1}=[\frac{t_j}{t_{j-1}}]$, the maximal integer not exceeding $\frac{t_j}{t_{j-1}}$,
$\sigma_{j,i}=[\frac{t_j-\sigma_{j,j-1}t_{j-1}-\cdots-\sigma_{j,i+1}t_{i+1}}{t_i}]$ for $2\le i \le j-2$,
and $\sigma_{j,1}=t_j-\sigma_{j,j-1}t_{j-1}-\cdots-\sigma_{j,2}t_2$.
Then $\sigma_{j,1}+\sigma_{j,2}t_2+\cdots+\sigma_{j,j-1}t_{j-1}=t_j$.
Let $\sigma_j=\sigma_{j,1}\ell_1+\sigma_{j,2}\ell_{t_2}+\cdots +\sigma_{j,j-1}\ell_{t_{j-1}}$ and let 
$$u_j=g^{\sigma_{j,1}}g_{\ell_{t_2}}^{\sigma_{j,2}}\cdots g_{\ell_{t_{j-1}}}^{\sigma_{j,j-1}}\in I_{\sigma_j,c}.$$
Since $\sigma_{j,j-1}\ge 1$, we have $\sigma_j>\ell_{t_{j-1}}$.
Further, since $\overline g_{\ell_{t_i}}=\overline c_{t_i}\overline g^{t_i}\in \overline A_c[\overline z]$ for every $i$ where $c_{t_i}\in A_c^*$ 
by Lemma \ref{Lemma 3.4}, $\overline u_j=\overline g_{\ell_{t_j}}\in \overline A_c[\overline z]$ up to units,
in particular, $\overline u_j=\overline g^{t_j}$ in $K[\overline z]$ up to units.
Hence it follows that $\ell_{t_{j-1}}<\sigma_j<\ell_{t_j}$ 
because $\overline I^K_{\ell_{t_j}}=\overline g^{t_j}K[\overline z]\supsetneq \overline I^K_{\ell_{t_j}+1}=\overline g^{t_j+1}K[\overline z]$. 
Since $\overline g_{\ell_{t_j}}=\overline d_j\overline u_j$ for $d_j \in A_c^*$, 
we have   
\begin{equation}
\alpha^{\ell_{t_j}}y_j =g_{\ell_{t_j}}=d_j u_j+\alpha \tilde f_j(z) \label{(10)}
\end{equation}  
where $\tilde f_j(z)\in I_{\sigma_j-1,c}$.  
Suppose $\ell_{t_j-i}<\sigma_j \leq \ell_{t_j-i+1}$ for some $1\le i < t_j-t_{j-1}$. 
Then by Lemma \ref{Lemma 3.5} (1), 
$I_{\sigma_j-1,c}=\alpha^{\sigma_j-1-\ell_{t_j-i}}I_{\ell_{t_j-i},c}+g_{\ell_{t_j-i+1}} A_c[z]$. 
Hence $\alpha \tilde f_j(z)\in \alpha^{\sigma_j-\ell_{t_j-i}}I_{\ell_{t_j-i},c}+\alpha g_{\ell_{t_j-i+1}} A_c[z]$.  
For $1\le l \le t_j-t_{j-1}$, since 
$$I_{\ell_{t_j-l},c}=\alpha^{\ell_{t_j-l}-\ell_{t_j-l-1}}I_{\ell_{t_j-l-1},c}+g_{\ell_{t_j-l}}A_c[z],$$
we have 
$$I_{\ell_{t_j-i},c}=\alpha^{\ell_{t_j-i}-\ell_{t_{j-1}-1}}I_{\ell_{t_{j-1}-1},c}+J_{t_j-i}$$ 
where $J_{t_j-i}=(\alpha^{\ell_{t_j-i}-\ell_{t_{j-1}}}g_{\ell_{t_{j-1}}}, \ldots,
\alpha^{\ell_{t_j-i}-\ell_{t_j-i-1}}g_{\ell_{t_j-i-1}}, g_{\ell_{t_j-i}})A_c[z]$. \\
Hence $g_{\ell_{t_j}}$ is written as 
\begin{align}
g_{\ell_{t_j}}= & d_j u_j+\alpha^{\sigma_j-\ell_{t_j-i}}g_{\ell_{t_j-i}}r_1+\alpha^{\sigma_j-\ell_{t_j-i-1}}g_{\ell_{t_j-i-1}}r_2  \nonumber \\
& + \cdots +\alpha^{\sigma_j-\ell_{t_{j-1}}}g_{\ell_{t_{j-1}}}r_{t_j-t_{j-1}}+\alpha^{\sigma_j-\ell_{t_{j-1}-1}}v_{\ell_{t_{j-1}-1}} 
+\alpha g_{\ell_{t_j-i+1}}\tilde r_j \label{(11)} 
\end{align}
where $r_1,\ldots, r_{t_j-t_{j-1}}, \tilde r_j\in A_c[z]$ and $v_{\ell_{t_{j-1}-1}}\in I_{\ell_{t_{j-1}-1},c}$. 
Note that $g_{\ell_{t_j-l}}\in I_{\ell_{t_j-l},c}$ for $1 \le l <t_j-t_{j-1}$ is written as 
$$g_{\ell_{t_j-l}}=g^{e_{l,1}}g_{\ell_{t_2}}^{e_{l,2}}\cdots g_{\ell_{t_{j-1}}}^{e_{l,j-1}}$$ 
where $e_{l,1}, \ldots, e_{l,j-1}$ are nonnegative integers such that
$e_{l,1}\ell_1+e_{l,2}\ell_{t_2}+\cdots +e_{l,j-1}\ell_{t_{j-1}}=\ell_{t_j-l}$ and $e_{l,1}+e_{l,2}t_2+\cdots+e_{l,j-1}t_{j-1}=t_j-l$.
Since
$$e_{l,1}+l+e_{l,2}t_2+\cdots+e_{l,j-1}t_{j-1}=t_j=\sigma_{j,1}+\sigma_{j,2}t_2+\cdots+\sigma_{j,j-1}t_{j-1},$$ 
it follows that $e_{l,j-1}\le \sigma_{j,j-1}$ by the definition of $\sigma_{j,j-1}$. 
With the notation above, the following assertion holds.

\begin{lem}\label{Lemma 3.7}
Let $j \ge 3$. Let $\tau_3=\sigma_3-\ell_{t_2}+(t_2-1)\ell_1$ and $\tau_j=\sigma_j-\ell_{t_{j-1}}+\tau_{j-1}$ for $j\ge 4$. 
Then 
$$g'_{\ell_{t_j}}=\alpha^{\tau_j}g' h_2 h_3\cdots h_j \quad \text{modulo} \; \alpha^{\tau_j+1}B_{\ell_{t_{j-1}},c}$$
where $h_j\in B_{\ell_{t_{j-1}},c}\setminus \alpha B_c$.
Hence
$$\delta(y_j)=\alpha^{q_j}\beta g'h_2\cdots h_j \quad \text{modulo} \; \alpha^{q_j+1}\beta B_{\ell_{t_{j-1}},c}$$
where $q_j=p-\ell_{t_j}+\tau_j\ge 0$. 
\end{lem}
\Proof
First, we consider the case $j=3$. Then $\ell_{t_2}< \sigma_3<\ell_{t_3}$.
Suppose $\ell_{t_3-i}<\sigma_3\le \ell_{t_3-i+1}$ for some $1 \le i < t_3-t_2$.  
By (\ref{(11)}), we have
\begin{align*}
g_{\ell_{t_3}}&=d_3 u_3+\alpha^{\sigma_3-\ell_{t_3-i}}g_{\ell_{t_3-i}}r_1+\alpha^{\sigma_3-\ell_{t_3-i-1}}g_{\ell_{t_3-i-1}}r_2 \\
& \quad + \cdots +\alpha^{\sigma_3-\ell_{t_2}}g_{\ell_{t_2}}r_{t_3-t_2}+\alpha^{\sigma_3-\ell_{t_2-1}}v_{\ell_{t_2-1}}+\alpha g_{\ell_{t_3}-i+1}\tilde r_3
\end{align*}
where $r_1,\ldots, r_{t_3-t_2}, \tilde r_3\in A_c[z]$ and $v_{\ell_{t_2-1}}\in I_{\ell_{t_2-1},c}$. \;
Since $g'_{\ell_{t_2}}=\alpha^{(t_2-1)\ell_1}\\ (g'h_2+\alpha b_1)$ for $b_1\in B_{\ell_1,c}$,  
\begin{align*}
u'_3 &=\sigma_{3,1}g^{\sigma_{3,1}-1}g_{\ell_{t_2}}^{\sigma_{3,2}}g'+\sigma_{3,2}g^{\sigma_{3,1}}g_{\ell_{t_2}}^{\sigma_{3,2}-1}g'_{\ell_{t_2}} \\
&=\alpha^{\sigma_3-\ell_1}\sigma_{3,1}y_1^{\sigma_{3,1}-1}y_2^{\sigma_{3,2}}g'+\alpha^{\sigma_3-\ell_{t_2}+(t_2-1)\ell_1}\sigma_{3,2}y_1^{\sigma_{3,1}}y_2^{\sigma_{3,2}-1}(g'h_2+\alpha b_1) \\
&=\alpha^{\tau_3}\sigma_{3,2}y_1^{\sigma_{3,1}}y_2^{\sigma_{3,2}-1}g'h_2 \quad \text{modulo}\; \alpha^{\tau_3+1}B_{\ell_{t_2},c} 
\end{align*}
and for $1 \le l<t_3-t_2$ 
\begin{align*}
g'_{\ell_{t_3-l}}=& e_{l,1}g^{e_{l,1}-1}g_{\ell_{t_2}}^{e_{l,2}}g'+e_{l,2}g^{e_{l,1}}g_{\ell_{t_2}}^{e_{l,2}-1}g'_{\ell_{t_2}} \\
=& \alpha^{\ell_{t_3-l}-\ell_1}e_{l,1}y_1^{e_{l,1}-1}y_2^{e_{l,2}}g' \\
& \hspace{2cm} +\alpha^{\ell_{t_3-l}-\ell_{t_2}+(t_2-1)\ell_1}e_{l,2}y_1^{e_{l,1}}y_2^{e_{l,2}-1}(g'h_2+\alpha b_1) \\
=& \alpha^{\ell_{t_3-l}-\ell_{t_2}+(t_2-1)\ell_1}e_{l,2}y_1^{e_{l,1}}y_2^{e_{l,2}-1}g'h_2 \\
& \hspace{4cm} \text{modulo}\; \alpha^{\ell_{t_3-l}-\ell_{t_2}+(t_2-1)\ell_1+1}B_{\ell_{t_2},c}.
\end{align*}
Hence we have  
\begin{align}
g'_{\ell_{t_3}}&=d_3 u'_3+\alpha^{\sigma_3-\ell_{t_3-i}}g'_{\ell_{t_3-i}}r_1+\cdots +\alpha^{\sigma_3-\ell_{t_2}}g'_{\ell_{t_2}}r_{t_3-t_2} \nonumber \\
&\qquad \quad +\alpha^{\sigma_3-\ell_{t_3-i}}g_{\ell_{t_3-i}}r'_1 +\cdots +\alpha^{\sigma_3-\ell_{t_2}}g_{\ell_{t_2}}r'_{t_3-t_2} \nonumber \\ 
& \qquad \quad +\alpha^{\sigma_3-\ell_{t_2-1}}v'_{\ell_{t_2-1}}+\alpha (g'_{\ell_{t_3}-i+1}\tilde r_3+g_{\ell_{t_3}-i+1}\tilde r'_3) \nonumber \\
&=\alpha^{\tau_3}g'h_2h_3 \quad \text{modulo}\;\alpha^{\tau_3+1}B_{\ell_{t_2},c} \label{(12)}
\end{align}
for 
$$h_3=d_3\sigma_{3,2}y_1^{\sigma_{3,1}}y_2^{\sigma_{3,2}-1}+b_3\in B_{\ell_{t_2},c}=A_c[z,y_1,y_2]$$
where $b_3\in A_c[z,y_1,y_2]$ satisfies $\alpha^{\sigma_3-\ell_{t_2}}b_3\in \alpha A_c[z]$ and $\deg_{y_2}h_3\le \sigma_{3,2}-1$.
Note that $e_{l,2}\le \sigma_{3,2}$. Note also that $v'_{\ell_{t_2-1}}\in \alpha^{\ell_{t_2-1}-\ell_1}B_{\ell_1,c}$. 
In fact, since $v_{\ell_{t_2-1}}=\alpha^{\ell_{t_2-1}}a(z,y_1)$ for $a(z,y_1)\in B_{\ell_1,c}$, we have
$$v'_{\ell_{t_2-1}}\delta(z)=\delta(v_{\ell_{t_2-1}})=\alpha^{\ell_{t_2-1}}(\partial_za \cdot \delta(z)+\partial_{y_1}a \cdot \delta(y_1)).$$
Hence it follows that $v'_{\ell_{t_2-1}}\in \alpha^{\ell_{t_2-1}-\ell_1}B_{\ell_1,c}$
since $\delta(z)=\alpha^p\beta$ and $\delta(y_1)=\alpha^{p-\ell_1}\beta g'$.

We show $h_3\notin \alpha B_c$. Suppose the contrary. 
Then $\tilde h_3:=\alpha^{\sigma_3-\ell_{t_2}}h_3\in I_{\sigma_3-\ell_{t_2}+1,c}$.
While, $\tilde h_3=d_3\sigma_{3,2}g^{\sigma_{3,1}}g_{\ell_{t_2}}^{\sigma_{3,2}-1}$ modulo $\alpha A_c[z]$, hence
$\overline{\tilde h_3}=\overline g^{t_3-t_2}$ in $K[\overline z]$ up to units. Thus it follows that $\sigma_3-\ell_{t_2}+1\le \ell_{t_3-t_2}$.
Write $\ell_{t_3-t_2}=e_1\ell_1+e_2\ell_{t_2}$ with nonnegative integers $e_1$ and $e_2$ such that $e_1+e_2t_2=t_3-t_2$.
Note that this expression is valid when $\ell_{t_3-t_2}=\ell_{t_2}$ or $\ell_1$. 
Then since $e_2+1\le \sigma_{3,2}$, we have
\begin{align*}
\sigma_3-\ell_{t_2}+1-\ell_{t_3-t_2} &=\sigma_{3,1}\ell_1+\sigma_{3,2}\ell_{t_2}-\ell_{t_2}+1-(e_1\ell_1+e_2\ell_{t_2}) \\
& =(t_3-\sigma_{3,2}t_2)\ell_1+\sigma_{3,2}\ell_{t_2}-\ell_{t_2}+1 \\
& \qquad \qquad \qquad -(t_3-t_2-e_2t_2)\ell_1-e_2\ell_{t_2} \\
& =(\sigma_{3,2}-e_2-1)(\ell_{t_2}-t_2\ell_1)+1>0,
\end{align*}
which is a contradiction.
Hence $h_3\notin \alpha B_c$.

By (\ref{(12)}), we have  
\begin{align*}
\alpha^{\ell_{t_3}}\delta(y_3) &=g'_{\ell_{t_3}}\delta(z)=\alpha^p \beta g'_{\ell_{t_3}} \\
&=\alpha^{p+\tau_3}\beta (g'h_2 h_3+\alpha \xi_3)
\end{align*}
for $\xi_3\in B_{\ell_{t_2},c}$.
Since $\beta g'h_2h_3\notin \alpha B_c$, it follows that $p+\tau_3\ge \ell_{t_3}$, i.e., $q_3\ge 0$ and
$$\delta(y_3)=\alpha^{q_3}\beta g'h_2h_3 \quad \text{modulo} \quad \alpha^{q_3+1}\beta B_{\ell_{t_2},c}.$$

Next, consider the case $j\ge 4$. We show by induction on $j$.
Suppose that $\ell_{t_j-i}<\sigma_j \leq \ell_{t_j-i+1}$ for some $1\le i < t_j-t_{j-1}$. 
Since $g'_{\ell_{t_s}}=\alpha^{\tau_s}g'h_2\cdots h_s$ modulo $\alpha^{\tau_s+1}B_c$ for $3\le s\le j-1$ and
$\sigma_j-\ell_{t_s}+\tau_s > \sigma_j-\ell_{t_{s+1}}+\tau_{s+1}$, we have  
\begin{align*}
u'_j =& \alpha^{\sigma_j-\ell_1}\sigma_{j,1}y_1^{\sigma_{j,1}-1}y_2^{\sigma_{j,2}}\cdots y_{j-1}^{\sigma_{j,j-1}}g' \\
& \qquad +\alpha^{\sigma_j-\ell_{t_2}}\sigma_{j,2}y_1^{\sigma_{j,1}}y_2^{\sigma_{j,2}-1}y_3^{\sigma_{j,3}}\cdots y_{j-1}^{\sigma_{j,j-1}}g'_{\ell_{t_2}}\\
& \qquad +\cdots +\alpha^{\sigma_j-\ell_{t_{j-1}}}\sigma_{j,j-1}y_1^{\sigma_{j,1}}\cdots y_{j-2}^{\sigma_{j,j-2}}y_{j-1}^{\sigma_{j,j-1}-1}g'_{\ell_{t_{j-1}}} \\
=&\alpha^{\tau_j}\sigma_{j,j-1}y_1^{\sigma_{j,1}}\cdots y_{j-2}^{\sigma_{j,j-2}}y_{j-1}^{\sigma_{j,j-1}-1}g'h_2\cdots h_{j-1}
\qquad \text{modulo} \quad \alpha^{\tau_j+1}B_{\ell_{t_{j-1}},c} 
\end{align*}
and for $1\le l <t_j-t_{j-1}$
\begin{align*}
g'_{\ell_{t_j-l}}=& \alpha^{\ell_{t_j-l}-\ell_{t_{j-1}}+\tau_{j-1}}e_{l,j-1}y_1^{e_{l,1}}\cdots y_{j-2}^{e_{l,j-2}}y_{j-1}^{e_{l,j-1}-1}
g'h_2\cdots h_{j-1} \\
& \hspace{5cm} \text{modulo} \quad \alpha^{\ell_{t_j-l}-\ell_{t_{j-1}}+\tau_{j-1}+1}B_{\ell_{t_{j-1}},c}.
\end{align*}  
By the induction hypothesis, 
\begin{equation}
\delta(y_s)=\alpha^{q_s}\beta g'h_2\cdots h_s \quad \text{modulo} \; \alpha^{q_s+1}B_{\ell_{t_{s-1}},c} \label{(13)}
\end{equation}  
for $3 \le s \le j-1$. 
Note that $q_s >q_{s+1}$. In fact, 
\begin{align*}
q_{s+1} &=p-\ell_{t_{s+1}}+\tau_{s+1} \\
& =p-\ell_{t_{s+1}}+\sigma_{s+1}-\ell_{t_s}+\tau_s \\
&=p-\ell_{t_s}+\tau_s-(\ell_{t_{s+1}}-\sigma_{s+1}) \\
&=q_s-(\ell_{t_{s+1}}-\sigma_{s+1})< q_s.
\end{align*}
By the same argument as in $j=3$,
\begin{align}
g'_{\ell_{t_j}}= & d_ju'_j +\alpha^{\sigma_j-\ell_{t_j-i}}g'_{\ell_{t_j-i}}r_1 + \cdots
+\alpha^{\sigma_j-\ell_{t_{j-1}}}g'_{\ell_{t_{j-1}}}r_{t_j-t_{j-1}} \nonumber \\
&\qquad +\alpha^{\sigma_j-\ell_{t_j-i}}g_{\ell_{t_j-i}}r'_1 +\cdots +\alpha^{\sigma_j-\ell_{t_{j-1}}}g_{\ell_{t_{j-1}}}r'_{t_j-t_{j-1}} \nonumber \\ 
& \qquad +\alpha^{\sigma_j-\ell_{t_{j-1}-1}}v'_{\ell_{t_{j-1}-1}}+\alpha (g'_{\ell_{t_j-i+1}}\tilde r_j+g_{\ell_{t_j-i+1}}\tilde r'_j) \nonumber \\
= & \alpha^{\tau_j}g'h_2\cdots h_j \quad \text{modulo}\;\alpha^{\tau_j+1}B_{\ell_{t_{j-1}},c} \label{(14)}
\end{align}
for 
$$h_j=d_j\sigma_{j,j-1}y_1^{\sigma_{j,1}}\cdots y_{j-2}^{\sigma_{j,j-2}}y_{j-1}^{\sigma_{j,j-1}-1}+b_j
\in B_{\ell_{t_{j-1}},c}=A_c[z,y_1,\ldots,y_{j-1}]$$
where $b_j\in A_c[z,y_1,\ldots,y_{j-1}]$ satisfies $\alpha^{\sigma_j-\ell_{t_{j-1}}}b_j\in \alpha A_c[z]$ and $\deg_{y_{j-1}}h_j\le \sigma_{j,j-1}-1$.
Note that $v'_{\ell_{t_{j-1}-1}}\in \alpha^{\ell_{t_{j-1}-1}-\ell_{t_{j-2}}+\tau_{j-2}}B_{\ell_{t_{j-2}},c}$ by (\ref{(13)}). 

Suppose $h_j \in \alpha B_c$. Then $\tilde h_j=\alpha^{\sigma_j-\ell_{t_{j-1}}}h_j\in I_{\sigma_j-\ell_{t_{j-1}}+1,c}$.
While, $\tilde h_j=d_j\sigma_{j,j-1}g^{\sigma_{j,1}}\cdots g_{\ell_{t_{j-2}}}^{\sigma_{j,j-2}}g_{\ell_{t_{j-1}}}^{\sigma_{j,j-1}-1}$
modulo $\alpha A_c[z]$, hence $\overline{\tilde h_j}=\overline g^{t_j-t_{j-1}}$ in $K[\overline z]$ up to units.
Thus we have $\sigma_j-\ell_{t_{j-1}}+1\le \ell_{t_j-t_{j-1}}$.
Write $\ell_{t_j-t_{j-1}}=e_1\ell_1+\cdots +e_{j-1}\ell_{t_{j-1}}$ with $e_1+\cdots +e_{j-1}t_{j-1}=t_j-t_{j-1}$.
Since $e_{j-1}+1\le \sigma_{j,j-1}$, it follows that 
\begin{align*}
&\sigma_j-\ell_{t_{j-1}}+1-\ell_{t_j-t_{j-1}} \\
=&\sigma_{j,1}\ell_1+\cdots +\sigma_{j,j-1}\ell_{t_{j-1}}-\ell_{t_{j-1}}+1-(e_1\ell_1+\cdots +e_{j-1}\ell_{t_{j-1}}) \\
= &(t_j-\sigma_{j,2}t_2-\cdots-\sigma_{j,j-1}t_{j-1})\ell_1+\sigma_{j,2}\ell_{t_2}+\cdots +\sigma_{j,j-1}\ell_{t_{j-1}}-\ell_{t_{j-1}}+1 \\
& \qquad \quad -(t_j-t_{j-1}-e_2t_2-\cdots -e_{j-1}t_{j-1})\ell_1-e_2\ell_{t_2}-\cdots -e_{j-1}\ell_{t_{j-1}} \\
=& (\sigma_{j,j-1}-e_{j-1}-1)(\ell_{t_{j-1}}-t_{j-1}\ell_1)+(\sigma_{j,j-2}-e_{j-2})(\ell_{t_{j-2}}-t_{j-2}\ell_1) \\
& \quad +\cdots +(\sigma_{j,2}-e_2)(\ell_{t_2}-t_2\ell_1)+1>0,
\end{align*}
which is a contradiction. Hence $h_j\notin \alpha B_c$.

By (\ref{(14)}), we have
$$\alpha^{\ell_{t_j}}\delta(y_j)=\alpha^p \beta g'_{\ell_{t_j}}=\alpha^{p+\tau_j}\beta (g'h_2 \cdots h_j+\alpha \xi_j)$$
for $\xi_j \in B_{\ell_{t_{j-1}},c}$.
Since $\beta g'h_2 \cdots h_j\notin \alpha B_c$, it follows that $q_j=p+\tau_j-\ell_{t_j}\ge 0$ and we obtain the expression of $\delta(y_j)$. 
\qed

\medskip

Let $j\ge 3$. We have by (\ref{(10)}) and (\ref{(11)}) $\alpha^{\ell_{t_j}}y_j=\alpha^{\sigma_j}\tilde h_j(z, y_1, \ldots, y_{j-1})$ for 
$$\tilde h_j(z, y_1, \ldots, y_{j-1})=d_jy_1^{\sigma_{j,1}}\cdots y_{j-1}^{\sigma_{j,j-1}}+\tilde b_j +\alpha a_j$$
where $a_j, \tilde b_j\in B_{\ell_{t_{j-1}},c}=A_c[z,y_1,\ldots,y_{j-1}]$ and
$\tilde b_j$ is a sum of terms $r(z)y_1^{e_1}\\ \cdots y_{j-1}^{e_{j-1}}$ such that
$e_1+e_2t_2+\cdots +e_{j-1}t_{j-1}<\sigma_{j,1}+\sigma_{j,2}t_2+\cdots +\sigma_{j,j-1}t_{j-1}$. 
Let $\tilde q_j=\ell_{t_j}-\sigma_j$. Then $\tilde q_j>0$ and we have   
$$\alpha^{\tilde q_j}y_j=\tilde h_j(z, y_1, \ldots, y_{j-1})\in B_{\ell_{t_{j-1}},c}.$$

So far, we obtain a sequence of $G_a$-equivariant affine modifications
$$A_c[z]\subset B_{\ell_1,c}\subset B_{\ell_{t_2},c} \subset \cdots \subset B_{\ell_{t_\nu},c}=B^{(1)}_c=A_c[z,y_1, \ldots, y_\nu]$$ 
with relations
\begin{align}
\alpha^{\ell_1}y_1 &=g(z), \nonumber \\
\alpha^{\tilde q_2}y_2 &=\tilde h_2(z, y_1), \nonumber \\
\alpha^{\tilde q_3}y_3 &=\tilde h_3(z, y_1,y_2), \label{(15)} \\
\cdots & \cdots \nonumber \\
\alpha^{\tilde q_\nu}y_\nu &=\tilde h_\nu(z, y_1, \ldots, y_{\nu-1}) \nonumber 
\end{align}
and
\begin{align}
\delta(z) &=\alpha^p \beta \nonumber \\
\delta(y_1) &= \alpha^{q_1}\beta g', \nonumber \\
\delta(y_2) &= \alpha^{q_2}\beta g'h_2 \quad \text{modulo} \; \alpha^{q_2+1}\beta A_c[z,y_1], \nonumber \\
\delta(y_3) &=\alpha^{q_3}\beta g' h_2h_3 \quad \text{modulo} \; \alpha^{q_3+1}\beta A_c[z,y_1,y_2], \label{(16)} \\
\cdots & \cdots \nonumber \\
\delta(y_\nu) &=\alpha^{q_\nu}\beta g' h_2\cdots h_\nu \quad \text{modulo} \; \alpha^{q_\nu+1}\beta A_c[z,y_1, \ldots, y_{\nu-1}]\nonumber 
\end{align}
where $h_i\in A_c[z,y_1, \ldots,y_{i-1}]\setminus \alpha B_c$ for $2 \le i \le \nu$.
We have shown that $0\le q_{i+1}< q_i$ for $3 \le i\le \nu-1$ in the proof of Lemma \ref{Lemma 3.7}. 
Since $q_1=p-\ell_1$, $q_2=q_1-\tilde q_2=p-\ell_1-\ell_{t_2}+t_2\ell_1$ and
$q_3=p-\ell_{t_3}+\tau_3=p-\ell_{t_2}+(t_2-1)\ell_1-(\ell_{t_3}-\sigma_3)$, we have
\begin{equation}
p>q_1>q_2>q_3>\cdots >q_\nu\ge 0. \label{(17)}
\end{equation}

Recall that there exists a sequence of $G_a$-equivariant affine modifications 
$$A[z]\subset B^{(1)}\subset B^{(2)} \subset \cdots \subset B^{(s)}=B.$$
Let $f_1, \ldots, f_d$ be the generators of $I^{(s)}=A[z]\cap \alpha_1^{\mu_1}\cdots \alpha_s^{\mu_s}B$.
Write $f_i=\alpha_1^{\mu_1}\cdots \alpha_s^{\mu_s}x_i$ for $x_i\in B$. Then
$$B=B^{(s)}=A[z][\alpha_1^{-\mu_1}\cdots \alpha_s^{-\mu_s}I^{(s)}]=A[z,x_1,\ldots,x_d].$$
Since $f_i\in I^{(s)}\subset I^{(1)}=I_\mu$ for every $i$ where $\mu=\mu_1$, we have $\alpha_2^{\mu_2}\cdots \alpha_s^{\mu_s}x_i\in B^{(1)}$. 
Hence by (\ref{(16)}),  
\begin{align*}
\alpha_2^{\mu_2}\cdots \alpha_s^{\mu_s}\delta(x_i)& \in (\delta(B^{(1)}))B \\
&\subset (\delta(B^{(1)}_c))B_c \\
& =(\delta(z), \delta(y_1), \ldots, \delta(y_\nu))B_c \\
& \subset \alpha_1^{q_\nu}B_c. 
\end{align*}
Since $B$ is factorial, we have $\delta(x_i)\in \alpha_1^{q_\nu}B$ for every $i$. 
Then it follows that $q_\nu=0$ since $\delta$ is irreducible. 
Hence we obtain the following results.

\begin{thm}\label{Theorem 3.8}
There exists $c \in A\setminus \alpha A$ so that $B^{(1)}_c=A_c[z, y_1, \ldots, y_\nu]$ with relations (\ref{(15)}) and (\ref{(16)}) 
where $q_\nu=0$.  
\end{thm}

\begin{thm}\label{Theorem 3.9}
Let $B$ be a factorial affine $k$-domain with an irreducible lnd $\delta$ and let $A=\Ker \delta$ be noetherian.
Suppose that $\delta(B)\cap A=\alpha^p\beta A$ where $p>0$ and $\alpha \in A$ is a prime element such that $\alpha$ and $\beta\in A$ are coprime. 
If $A$ is an affine $k$-domain, then the general closed fiber of $\pi : X=\Spec B \to Y=\Spec A$ over $\Gamma =\Spec A/\alpha A$ 
consists of a disjoint union of $m$ affine lines where $m\ge \deg_{\overline z}\overline g(\overline z)$. 
\end{thm}
\Proof
Since $\delta$ is irreducible, the induced lnd $\overline \delta$ on $\overline B=B/\alpha B$ is nontrivial. 
The restriction $\pi |_{\Gamma} : \Spec \overline B\to \Gamma$ is dominant by \cite[Lemma 1.6]{GMM} and 
decomposes as $\pi |_{\Gamma}=\tau \circ \rho$ where $\rho : \Spec \overline B\to \Spec (\Ker \overline \delta)$ is the quotient morphism and
$\tau : \Spec (\Ker \overline \delta)\to \Spec \overline A$ is the morphism induced by the inclusion
$\overline A \hookrightarrow \Ker \overline \delta$. 
By Theorem \ref{Theorem 3.8}, we have $B_{c\beta}=A_{c\beta}[z, y_1, \ldots, y_\nu]$ for some $c\in A\setminus \alpha A$ with relations 
(\ref{(15)}) and (\ref{(16)}) where $q_\nu=0$. 
Hence $\overline B_{\overline{c\beta}}=\overline A_{\overline{c\beta}}[\overline z, \overline y_1, \ldots, \overline y_\nu]$ and
$\overline y_\nu$ is a local slice of $\overline \delta$ by (\ref{(16)}) and (\ref{(17)}). 
Since the general closed fiber of $\rho$ is $\A^1$, it follows from (\ref{(15)}) that
the general closed fiber of $\pi |_{\Gamma}$ consists of disjoint union of $m$ affine lines where
$m=[Q(\Ker \overline \delta):Q(\overline A)]\ge \deg_{\overline z}\overline g(\overline z)$. 
\qed

\medskip

By Lemma \ref{Lemma 3.2}, $\deg_{\overline z}\overline g(\overline z)>1$.
Hence the general closed fiber of the quotient morphism $\pi : X=\Spec B \to Y=\Spec A$ over $\Gamma_i=\Spec A/\alpha_i A$
consists of a disjoint union of $m_i$ affine lines where $m_i\ge 2$. 

Suppose $p=1$. Then $\nu=1$ by (\ref{(17)}) and we have a $G_a$-equivariant affine modification 
$$A_c[z]\subset B^{(1)}_c=A_c[z,y_1]$$
with $\delta(z)=\alpha \beta$ and $\delta(y_1)=\beta g'(z)$. 

Consider the case $s=1$, i.e., $\delta(B)\cap A=\alpha^p A$. 
Then $B=B^{(1)}$ and $B_c=A_c[z,y_1, \ldots, y_\nu]$ with relations (\ref{(15)}) and (\ref{(16)}) where $\beta=1$. 
If $\ell_1=p$, i.e., $g(z)=\alpha^p y_1$ and $\overline A$ is factorial, then we have $q_1=0$ and $A[z]\subset B^{(1)}=A[z,y_1]=B$ (cf. \cite{Ma1}).
Hence $B\cong A[Y,Z]/(\alpha^p Y-g(Z))$ where $A[Y,Z]=A^{[2]}$. 
Let $e\ge 0$ be the minimal integer such that $\theta :=c^e g' h_2\cdots h_\nu \in B$. 

\begin{cor}\label{Corollary 3.10}
Suppose that $\delta(B)\cap A=\alpha^p A$ for a prime element $\alpha \in A$ and $p>0$. 
With the notation above, the following assertions hold. 
\begin{enumerate}
\item[(1)] If $\overline{\theta}$ is a unit of $\overline B=B/\alpha B$, 
then $X=\Spec B$ has no fixed points under the $G_a$-action corresponding to $\delta$. 
\item[(2)] Suppose that $\overline A$ is factorial. 
If $\ell_1=p$, then $X$ has no fixed points under the $G_a$-action corresponding to $\delta$ if and only if $\overline{g'}\in {\overline B}^*$. 
In particular, if $p=1$ and $\overline{g'}\in {\overline B}^*$, then $X^{G_a}=\emptyset$.
Further, if ${\overline B}^*={\overline A}^*$, then $X^{G_a}\neq \emptyset$. 
\end{enumerate}
\end{cor}
\Proof
(1) It follows that $\sqrt{(\delta (B))B}\supset \sqrt{(\delta(z), c^e\delta (y_\nu))B} \supset (\alpha, \theta)B$. 
Hence the fixed point locus $X^{G_a}$ is contained in the closed set $V(\alpha, \theta)$. 
Since $\overline{\theta}\in {\overline B}^*$, $(\alpha, \theta)B$ is a unit ideal and $V(\alpha, \theta)=\emptyset$.
Thus $X^{G_a}=\emptyset$ follows. 

(2) Since $B=A[z,y_1]$ with $\delta(z)=\alpha^p$ and $\delta(y_1)=g'(z)$, 
it follows that $(\delta(B))B=(\delta(z), \delta(y_1))B=(\alpha^p, g')B$. 
Hence the first assertion follows. 

Suppose ${\overline B}^*={\overline A}^*$. Then $\overline{g'}\in \overline B^*$ implies $\overline{g'}=\overline a$ for $a\in A\setminus \alpha A$.
Then we have $g'-a \in I_1$, which is a contradiction by Lemma \ref{Lemma 3.2}. 
Hence $\overline{g'}\notin {\overline B}^*$ and $X^{G_a}\neq \emptyset$. 
\qed

\medskip

\Remark Note that $\overline{g'} \in \overline B$ is a factor of $\overline a \in \overline B$ for some $a \in A\setminus \alpha A$.
In fact, since $\overline g\in \overline A[\overline z]\subset K[\overline z]$ is irreducible,
we have $\overline g \eta_1 +\overline g' \eta_2=1$ for $\eta_1, \eta_2 \in K[\overline z]$.
Hence $\overline g \zeta_1 +\overline g' \zeta_2=\overline a$ in $\overline A[\overline z]$ for $a\in A \setminus \alpha A$
and $\zeta_1, \zeta_2 \in \overline A[\overline z]$. Since $\overline g=0$ in $\overline B$, we have $\overline g' \zeta_2=\overline a$. 

\bigskip

Note that if $\overline B^*=k^*$, then $\overline B^*=\overline A^*=k^*$. 

\begin{cor}\label{Corollary 3.11}
Let $X=\Spec B$ be a smooth factorial affine variety with a $G_a$-action associated to an irreducible lnd $\delta$. 
Suppose that $A=\Ker \delta$ is an affine $k$-domain and $Y=\Spec A$ is smooth. 
Suppose further that the quotient morphism $\pi : X \to Y$ is surjective and equi-dimensional, and
the restriction $\pi |_{\pi^{-1}(D(a))} : \pi^{-1}(D(a)) \to D(a)$ is a trivial $\A^1$-bundle for a nonzero $a\in A$. 
If the general closed fiber of $\pi$ over $V(a)$ is irreducible, then $\pi : X \to Y$ is a trivial $\A^1$-bundle.
Hence $X$ is equivariantly isomorphic to $Y \times \A^1$ 
where $G_a$ acts trivially on $Y$ and by translation on $\A^1$.  
\end{cor}
\Proof
By Lemma \ref{Lemma 2.1}, $\delta(B)\cap A$ is a principal ideal.
If $\delta(B)\cap A$ is a unit ideal, then there exsits a slice of $\delta$ and $B=A^{[1]}$, i.e., $\pi$ is a trivial $\A^1$-bundle. 
Suppose that $\delta(B)\cap A=a'A$ for a non-unit $a' \in A$. 
Let $a_1\in A$ be any prime factor of $a'$. 
If $V(a_1)\not\subset V(a)$, then the general closed fiber over $V(a_1)$ consists of a single $\A^1$
since $\pi|_{\pi^{-1}(D(a))}$ is a trivial $\A^1$-bundle. 
However, this is a conrtadiction by Theorem \ref{Theorem 3.9}.
Hence $V(a_1)\subset V(a)$, and it follows that $V(a')\subset V(a)$. 
Since the general closed fiber of $\pi$ over $V(a)$ is irreducible by the assumption,
this is a contradiction again by Theorem \ref{Theorem 3.9}. 
Hence the assertion follows.
\qed

\medskip

We apply the results obrained so far to an affine pseudo-$n$-space. 
Let $n\ge 3$. An affine pseudo-$n$-space $X=\Spec B$ is a smooth affine algebraic variety equipped with a faithfully flat morphism
$q : X \to \A^1=\Spec k[x]$ such that $q^{-1}(\A^1_*)\cong \A^1_*\times \A^{n-1}$ and $q^*(0)$ is irreducible and reduced.  
Then $B$ is factorial, $B^*=k^*$, and $x$ is a prime element of $B$ by \cite{Ma1}. 
The following can be proved by the same argument in \cite[Theorem 2.3]{GMM2}.

\begin{thm}\rm{(\cite[Theorem 2.3]{GMM2})}\label{Theorem 3.12} 
Let $X=\Spec B$ be an affine pseudo-$n$-space with a faithfully flat morphism $q : X \to \A^1=\Spec k[x]$.
Assume that $X_0=q^*(0)$ is smooth.
Then the following conditions are equivalent.
\begin{enumerate}
\item[(1)] $X$ is contractible.
\item[(2)] $X$ is acyclic, i.e., $H_i(X;\Z)=0$ for every $i>0$.
\item[(3)] $X_0$ is acyclic.
\end{enumerate}
\end{thm}

\medskip
  
By the trivialization $q^{-1}(\A^1_*)\cong \A^1_*\times \A^{n-1}$,
$X$ has an algebraic action of $G_a^{n-1}$ associated to commuting irreducible lnds $\delta_1, \ldots, \delta_{n-1}$ such that
$\cap_{i=1}^{n-1}\Ker \delta_i=k[x]$. 
For each $i$, there exsits a local slice $z_i\in B$ such that $\delta_i(z_i)=x^{p_i}$ for $p_i\ge 0$ \cite{Ma1}.
Take any lnd $\delta_i$ and let $\delta=\delta_i$.
Then $X$ has a $G_a$-action corresponding to $\delta$ and $q : X \to \A^1$ is $G_a$-equivariant. 
Suppose that $A=\Ker \delta$ is an affine $k$-domain. Then $q$ splits to $q=\tau\circ \pi$ where $\pi : X \to Y=\Spec A$ is the quotient morphism
and $\tau : Y \to \A^1=\Spec k[x]$ is the morphism induced by the inclusion $k[x]\hookrightarrow A$. 
If $Y$ is smooth and $\pi$ is surjective and equi-dimensional, $\delta(B)\cap A$ is principal by Lemma \ref{Lemma 2.1} and 
$\delta(B)\cap A=x^p A$ for $p\ge 0$. 
Applying Corollary \ref{Corollary 3.11} to an affine pseudo-$n$-space $X$, we have the following. 

\begin{cor}\label{Corollary 3.13}
Let $X=\Spec B$ be an affine pseudo-$n$-space with a faithfully flat morphism $q : X \to \A^1=\Spec k[x]$ and a $G_a$-action 
such that $q: X \to \A^1$ is $G_a$-equivariant where $G_a$ acts trivially on $\A^1$. 
Let $\delta$ be the irreducible lnd on $B$ corresponding to the $G_a$-action and let $A=\Ker \delta$.
Suppose that $A$ is an affine $k$-domain and $Y=\Spec A$ is smooth.
Suppose further that the quotient morphism $\pi : X \to Y$ is surjective and equi-dimensional. 
If the general closed fiber of $\pi$ over $V(x)$ is irreducible, 
then $X$ is $G_a$-equivariantly isomorphic to $Y \times \A^1$ where $G_a$ acts on $\A^1$ by translation.  
\end{cor}

\medskip

For an affine pseudo-$3$-space $X$, the following holds.

\begin{cor}\label{Corollary 3.14}
Let $X=\Spec B$ be as in Corollary \ref{Corollary 3.13} with $n=3$. 
Suppose that $q^*(0)$ is factorial with $(B/x B)^*=k^*$.
Then $X\cong \A^3$ and $x$ is a variable. 
Further, if the general closed fiber of $\pi$ over $V(x)$ is irreducible, then $X\cong Y\times \A^1$ where $Y=\Spec A$ with $A=B^{G_a}=k[x]^{[1]}$. 
\end{cor}
\Proof
As shown above, there exist an irreducible lnd $\delta$ and $z\in B$ such that $\delta(z)=x^{p_1}$ for $p_1\ge 0$.
The kernel $A=\Ker \delta$ is $k[x]^{[1]}$ by \cite[Lemma 5.1]{Ma1}. 
Since $q^*(0)=\Spec B/xB$ is factorial with a nontrivial $G_a$-action induced by $\delta$ and $(B/xB)^*=k^*$, it follows that $q^*(0)\cong \A^2$,
and hence $q^*(0)$ is smooth and acyclic.
Then $X$ is acyclic by Theorem \ref{Theorem 3.12}. 
Since $q^*(0)$ is factorial, $X \cong \A^3$ and $x$ is a variable by a result of Kaliman \cite{Kal2}. 
Also by \cite{Kal1}, $\pi : X \to Y=\Spec A\cong \A^2$ is surjective and equi-dimensional. 
Hence the plinth ideal $\delta(B)\cap A$ is principal by Lemma \ref{Lemma 2.1} and generated by $x^p$ for $p\ge 0$.
If the general closed fiber of $\pi$ over $V(x)$ is irreducible, $X\cong Y\times \A^1$ by Corollary \ref{Corollary 3.13}.
\qed

\medskip

In Corollary \ref{Corollary 3.14}, if the general closed fiber of $\pi$ over $V(x)$ is reducible, 
then $X$ is not necessarily isomorphic to $Y\times \A^1$ as shown in Example 4.1 below, although $X\cong \A^3$.

\medskip
\bigskip

\section{Examples}

We illustrate affine modifications of $G_a$-varieties by examples.

\medskip

\noindent{\bf Example 4.1}\label{Example 4.1}

Let $B=k[x,y,z]$ be a polynomial ring with an lnd $\delta$ defined by $\delta(x)=0$, $\delta(y)=-2z$, $\delta(z)=x^2$.
Then $A=\Ker \delta=k[x,t]=k^{[2]}$ where $t=x^2y+z^2$ and $\delta(B)\cap A=x^2 A$. 
It follows from $\delta(z)=x^2$ that $B[x^{-1}]=A[x^{-1}][z]$. 
Let $g(z)=z^2-t$. Then $g(z)=-x^2y$ and $I_1=A[z]\cap xB=(x,g)A[z]$. 
Since $I_2=A[z]\cap x^2B=(x^2,g)A[z]$, we have $B_2=A[z][x^{-2}I_2]=B$, and a sequence of $G_a$-equivariant affine modifications
$$A[z]\subset B_1=A[z][x^{-1}I_1]=A[z, xy]\subset B_2=B.$$
Let $\pi : X=\Spec B \to Y=\Spec A$ be the quotient morphism. 
The singular locus $\Sing (\pi)=\{Q\in Y \mid \pi^*(Q)\not\cong \A^1_{k(Q)}\}$ consists of one irreducible component
$V(x)\cong \A^1$ of $Y=\A^2_{(x,t)}$ where $k(Q)$ is the residue field of $Y$ at $Q$. 
For $Q=(0, \beta)\in V(x)$, the fiber $\pi^*(Q)$ is $\A^1+\A^1$ if $\beta\neq 0$ and $2\A^1$ if $\beta=0$. 
The fixed point locus consists of the fiber $\pi^{-1}(O)$ for $O=(0,0)\in Y$.
Note that the degree of $g(z)\in A[z]$ modulo $xA[z]$ is two. 
Though $\pi : X=\A^3 \to Y=\A^2$ is not a trivial $\A^1$-bundle, $x$ is a variable of $B$. 

\medskip

\noindent{\bf Example 4.2}\label{Example 4.2}

Let $m\ge 0$ and let $B=R[x,y,z]$ be a polynomial ring over $R=k^{[m]}$ with a locally nilpotent $R$-derivation $\delta$ defined by
$$\delta(x)=0, \quad \delta(y)=h(x,z), \quad \delta(z)=f(x)^p$$
where $p>0$, $f(x)\in R[x]\setminus R$ and $h(x,z)\in R[x,z]\setminus R[x]$.
We assume that $f(x)$ is irreducible and $(f(x), a(x))B$ is a unit ideal
where $a(x)\in R[z]$ is the coefficient of the highest term of $h(x,z)$ with respect to $z$. 
Then the lnd $\delta$ is irreducible and its kernel $A=\Ker \delta$ is $R[x, F]=k^{[m+2]}$ 
where $F=f(x)^p y-g(x,z)$ and $g(x,z)\in R[x,z]$ is a polynomial such that $\partial_z g(x,z)=h(x,z)$. 
If $\delta(B)\cap A$ is a unit ideal, there exists a slice $s\in B$ and 
$B=A[s]=R[x, F, s]$. Hence $F=f(x)^p y-g(x,z)$ is an $x$-variable of $B=R[x,y,z]$.
We assume $(f(x), h(x,z))B\neq B$. Then the fixed point locus of $X=\Spec B$ under the $G_a$-action corresponding to $\delta$ is nonempty and
$\delta(B)\cap A\neq A$. 
The quotient morphism $\pi : \A^{m+3}=\Spec B \to Y=\Spec A=\A^{m+2}$ is surjective and equi-dimensional. 
Hence the plinth ideal $\delta(B)\cap A$ is principal and $\delta(B)\cap A=f^p A$. 
In fact, the plinth ideal $\delta(B)\cap A$ is generated by $f^{p'}$ for $0<p' \le p$ since $\delta(z)=f^p$. 
Suppose $p'<p$. Then there exsits $\xi(x,y,z)\in B$ such that 
$\delta(\xi)=\partial_y \xi \cdot \delta(y)+\partial_z \xi \cdot \delta(z)=\partial_y\xi \cdot h(x,z)+\partial_z\xi \cdot f^p=f^{p'}$.
Since $h(x,z)$ is not divisible by $f$, it follows that $(f^{p-p'}, h)B$ is a unit ideal, which contradicts to the assumption $(f,h)B\neq B$. 
Hence $\delta(B)\cap A=f^p A$. 
Then it follows that $I_1=A[z]\cap f B=(f,g+F)A[z]$ and $I_p=A[z]\cap f^p B=(f^p, g+F)A[z]$. 
Hence we have $B_p=A[z][f^{-p}I_p]=A[z,y]=B$ and a sequence of $G_a$-equivariant affine modifications 
$$A[z]=R[x,f^py,z]\subset B_p=B=R[x,y,z].$$ 
The closure of $\Sing (\pi)$ consists of the single irreducible component $V(f)$ and
the general closed fiber of $\pi$ over $V(f)$ consists of $m$ affine lines where $m$ is the degree of $g+F\in A[z]$ modulo $fA[z]$ 
which coincides with $\deg_z h(x,z)+1$.

\medskip

\noindent{\bf Example 4.3}\label{Example 4.3} 

Let $B=k[x,y,z]=k^{[3]}$. Consider the lnd $\delta$ on $B$ defined by 
$$\delta(x)=-2FR, \quad \delta(y)=6x^2R-G, \quad \delta(z)=2x(5yR+F^2)$$
where 
$$F=xz-y^2, \quad G=zF^2+2x^2yF+x^5, \quad R=x^3+yF.$$
This lnd was studied by Freudenburg and is called the $(2,5)$ derivation 
(\cite{Fre}, \cite{Fre1}). 
The lnd $\delta$ is irreducible and its kernel $A=\Ker \delta$ is $k[F,G]=k^{[2]}$. 
The plinth ideal $\delta(B)\cap A$ is generated by $FG$, and $\delta(R)=-FG$ \cite{Fre1}. 
Though $A/F A\cong A/G A \cong k^{[1]}$ is factorial, neither $B/FB$ nor $B/GB$ is factorial. 
The fixed point locus $X^{G_a}$ of $X=\Spec B$ is nonempty and defined by $x=y=0$.
Note that there are relations 
$$R^2+F^3=Gx, \quad FS=G-x^2 R \quad \text{where} \quad S=x^2y+Fz.$$ 
Hence $Gx$, $FG^2S$, $FG^3y=G^3(R-x^3)$ and $F^2G^5z=FG^5(S-x^2y)$ are elements of $A[R]=k[F,G,R]$.

Put $\alpha_1=F$ and $\alpha_2=G$ with the notation in section 3. Then $I^{(1)}=A[R]\cap F^2B$ and 
$$I_1=A[R]\cap F B=(F, g)A[R]$$
where $g=FG^2S=G^3-(R^2+F^3)^2R$.
Note that $v=FG^3y=gR-(R^2+F^3)^2F^3\in I_1$ and $w=F^2G^5z=g^2+(R^2+F^3)^4F^3\in I^{(1)}\subset I_1$. 
We have  
$$B^{(1)}=B_1=A[R][F^{-1}I_1]=k[F,G,R,G^2S].$$
Note that $G^3y, G^5z\in B^{(1)}$.
We have 
$$I^{(2)}=A[R]\cap F^2G^5 B, \quad B^{(2)}=A[R][F^{-2}G^{-5}I^{(2)}]=B.$$
Note that $x, y, z\in B^{(2)}$ since $F^2G^4u, FG^2v, w \in I^{(2)}$ where $u=R^2+F^3=Gx$.
The sequence of $G_a$-equivariant affine modifications is 
$$A[R]=k[F,G,R]\subset B^{(1)}=k[F,G,R,G^2S]\subset B^{(2)}=B.$$ 

If we put $\alpha_1=G$ and $\alpha_2=F$, then $I^{(1)}=A[R]\cap G^5 B$ and 
$$I_1=A[R]\cap G B=(G, u)A[R].$$
Hence
$$B^{(1)}=A[R][G^{-1}I_1]=k[F, G, R, x].$$
Note that $FS=G-x^2R\in B^{(1)}$, $Fy=R-x^3 \in B^{(1)}$, and hence $F^2z=FS-x^2Fy \in B^{(1)}$. 
We have $I^{(2)}=A[R]\cap F^2G^5 B$ and $B^{(2)}=A[R][F^{-2}G^{-5}I^{(2)}]=B$. 
The sequence of $G_a$-equivariant affine modifications is 
$$A[R]=k[F,G,R]\subset B^{(1)}=k[F,G,R,x]\subset B^{(2)}=B.$$ 

Let $\pi : X \to Y=\Spec k[F,G]$ be the quotient morphism.
The singular locus $\Sing (\pi)$ consists of $2$ components $V(F)$ and $V(G)$.
The closed fiber of $\pi$ over $V(F)\setminus V(F,G)$ consists of $5$ affine lines, the one over $V(G)\setminus V(F,G)$ consists of $2$ affine lines, 
and the one over $V(F,G)$ is $\A^1$ with multiplicity $10$, which is the fixed-point locus. 
Note that the degree of $g\in A[R]$ modulo $FA[R]$ is five and the one of $u\in A[R]$ modulo $GA[R]$ is two. 

The canonical factorization of $X$ is explicitly given by Freudenburg in \cite{Fre1};
\begin{align*}
A=k[F,G]\subset A[R]=& k[F,G,R]\subset k[F,G,R,x] \\
  & \subset k[F,R,x,S] \subset k[F,x,S,y]\subset k[x,y,z]=B.
\end{align*}  
The subsequence $A[R]\subset \cdots \subset B$ of the canonical factorization consists of $G_a$-equivariant affine modifications.

Our sequences of $G_a$-equivariant affine modifications can provide the information on the singular fibers of the quotient morphism in a direct way.

\end{document}